%% file: gosnip.tex
\documentclass{article}

\usepackage{dpr_fec,hyperref,enumitem,graphics,subcaption,fullpage}
\usepackage{tikz}
\usetikzlibrary{shapes,arrows}
\pgfdeclarelayer{background}
\pgfdeclarelayer{foreground}
\pgfsetlayers{background,main,foreground}
\usepackage[titletoc]{appendix}
\newcommand{\background}[5]{%
  \begin{pgfonlayer}{background}
    \path (#1.west |- #2.north)+(-0.5,0.25) node (a1) {};
    \path (#3.east |- #4.south)+(+0.5,-0.25) node (a2) {};
    \path[rounded corners, draw=black!50, dashed]
      (a1) rectangle (a2);
      \path (#3.east |- #2.north)+(0,0.25)--(#1.west |- #2.north) node[midway] (#5-n) {};
      \path (#3.east |- #2.south)+(0,-0.35)--(#1.west |- #2.south) node[midway] (#5-s) {};
      \path (#3.east |- #2.north)+(0.7,0)--(#3.east |- #4.south) node[midway] (#5-w) {};
  \end{pgfonlayer}}
  
\hypersetup{colorlinks,urlcolor=blue}
\allowdisplaybreaks

\usepackage{natbib}
 \bibpunct[, ]{(}{)}{,}{a}{}{,}%
\begin{document}

\title{A Decomposition Algorithm for Large-Scale Security-Constrained AC Optimal Power Flow}

\author{%
Frank E. Curtis\thanks{Department of Industrial and Systems Engineering, Lehigh University, Bethlehem, PA, USA, \href{mailto:frank.e.curtis@gmail.com}{frank.e.curtis@gmail.com}} \and
Daniel K. Molzahn\thanks{School of Electrical and Computer Engineering, Georgia Institute of Technology, Atlanta, GA, USA, \href{mailto:molzahn@gatech.edu}{molzahn@gatech.edu}} \and
Shenyinying Tu\thanks{Department of Industrial Engineering and Management Sciences, Northwestern University, Evanston, IL, USA, \href{mailto:shenyinyingtu2021@u.northwestern.edu}{shenyinyingtu2021@u.northwestern.edu}} \and
Andreas W\"achter\thanks{Department of Industrial Engineering and Management Sciences, Northwestern University, Evanston, IL, USA, \href{mailto:waechter@iems.northwestern.edu}{waechter@iems.northwestern.edu}} \and
Ermin Wei\thanks{Department of Electrical and Computer Engineering, Department of Industrial Engineering and Management Sciences, Northwestern University, Evanston, IL, USA, \href{mailto:ermin.wei@northwestern.edu}{ermin.wei@northwestern.edu}} \and
Elizabeth Wong\thanks{Department of Mathematics, University of California, San Diego, La Jolla, CA, USA, \href{mailto:elwong@ucsd.edu}{elwong@ucsd.edu}}
}

\maketitle

\begin{abstract}
  \input{abstract}
\end{abstract}


\input{body}


\section*{Acknowledgment}

The authors would like to acknowledge Xingyi He for the training of the machine learning model for the initial ranking of contingencies. The authors also gratefully acknowledge the organizers of the ARPA-E Grid Optimization competition. The information, data, or work presented herein was funded in part by the Advanced Research Projects Agency-Energy (ARPA-E), U.S. Department of Energy, under Award Number DE-AR0001073.

\bibliographystyle{natbib}
\bibliography{gosnip}

\end{document}

%% file: abstract.tex
A decomposition algorithm for solving large-scale security-constrained AC optimal power flow problems is presented.  The formulation considered is the one used in the ARPA-E Grid Optimization (GO) Competition, Challenge~1, held from November 2018 through October 2019.  The techniques found to be most effective in terms of performance in the challenge are presented, including strategies for contingency selection, fast contingency evaluation, handling complementarity constraints, avoiding issues related to degeneracy, and exploiting parallelism.  The results of numerical experiments are provided to demonstrate the effectiveness of the proposed techniques as compared to alternative strategies.

%% file: body.tex
\newcommand{\ouralg}{\textbf{Algorithm GO-SNIP}}

\section{Introduction}\label{sec.introduction}

Operators of electric power systems solve optimal power flow (OPF) problems in order to manage the transfer of power from generators to loads through an electric grid efficiently and reliably. First formulated in~\citep{carpentier1962}, OPF problems seek minimum-cost operating points that satisfy both engineering constraints (e.g., line flow, voltage magnitude, and power generation limits) and the power flow equations, which model the physics of the transmission network. The most prominent use of OPF problems is scheduling generator setpoints by clearing electricity markets~\citep{kirschen2018fundamentals,ferc1}. Closely related optimization problems are also employed for a wide variety of other purposes, such as computing actions to restore acceptable operation during emergencies~\citep{misra2017}, managing uncertainty from stochastic generators~\citep{bienstock2014chance}, determining stability margins~\citep{avalos2009}, and planning system expansions~\citep{mahdavi2019}, amongst others.

OPF problems are challenging due to their size and mathematical characteristics. Industrial-scale OPF problems consider power systems involving thousands of generators and tens of thousands of buses and transmission lines, leading to commensurate numbers of variables and constraints. Regarding their mathematical characteristics, the steady-state physics of power systems are most accurately modeled using the nonlinear system of \emph{AC power flow} equations. AC-OPF problems, which incorporate AC power flow as equality constraints, are nonconvex, may have multiple local optima~\citep{bukhsh2013}, and are NP-hard \citep{vanhentenryck2016,bienstock2015}.

The need to maintain power system reliability imposes additional challenges, leading to so-called \emph{security-constrained} generalizations of OPF problems~\citep{stott1987}. Power systems are operated in order to withstand component failures within a specified set of \emph{contingency} scenarios, which typically include the individual failures of each component, i.e., \emph{$N-1$ security}~\citep{doe2016}. The large number of contingencies significantly increases the computational burden as each contingency ostensibly introduces a similar number of variables and constraints as the nominal problem. Moreover, for the security-constrained OPF formulation considered in this paper, the recourse model for the contingencies includes \emph{complementarity} conditions related to the limits on the power outputs for each generator. These complementarity conditions introduce additional nonlinearities to the formulation, as well as certain degeneracies.

Despite these challenges, operators require solutions to security-constrained AC-OPF (SC-AC-OPF) problems within demanding time constraints, e.g., approximately 10 minutes, for many real-time applications. Rather than attempt to solve SC-AC-OPF problems directly, operators often resort to a so-called \emph{DC power flow} approximation of AC power flow.  A DC power flow approximation applies various assumptions based on characteristics of typical transmission systems to linearize the power flow equations. While this provides advantages in the reliability and speed of the computations, a DC power flow approximation neglects important aspects of the power system physics, including voltage magnitudes and reactive power.  Hence, DC power flow is inapplicable for some uses of OPF problems. Moreover, a DC power flow approximation can introduce significant errors that impact both the achievable cost and the feasibility of the resulting DC-OPF solutions~\citep{stott2009,dvijotham_molzahn-cdc2016}.

Replacing existing industrial applications of DC-OPF problems with more accurate AC-OPF problems has the potential to reduce operating costs and constraint violations significantly, with plausible annual savings ranging from \$6~billion to \$9~billion in the United States alone~\citep{ferc1}. However, power systems engineers have often considered industrial-scale applications of AC-OPF problems to be intractable due to the aforementioned computational challenges. 

Over the last two decades, the research community has developed many new techniques for modeling and solving AC-OPF problems~\citep{capitanescu2011,capitanescu2016,molzahn_hiskens-fnt2019}. In order to evaluate and compare the progress in this area of research, the Advanced Research Projects Agency--Energy (ARPA-E) in the US Department of Energy organized the Grid Optimization (GO) Challenge 1 competition, which was held from November 2018 through October 2019, for solving SC-AC-OPF problems~\citep{GoCompetition}. Each team participating in this competition submitted an implementation of an SC-AC-OPF solution algorithm to the competition organizers.  The organizers ran each implementation for numerous SC-AC-OPF test cases using the same computing hardware and scored the outputs.

This paper presents the algorithm that our team, ``GO-SNIP,'' developed for the GO Challenge~1 competition along with some subsequent improvements.  Our algorithm uses an interior-point optimization algorithm as its core, and confronts large problem sizes by iteratively enforcing the constraints associated with contingencies that the algorithm identifies to be important.  Our algorithm applies several contingency screening and parallel processing techniques to identify quickly what seem to be the most important contingencies.  Additionally, our algorithm utilizes tailored heuristics designed to manage the complementarity conditions and avoids certain degeneracies by modifying the original problem formulation.

For the GO Challenge 1 competition, our algorithm successfully solved SC-AC-OPF problems for 17 synthetic networks (each with 20 scenarios) and 3 actual industrial networks (each with 4 scenarios) for a total of 352 SC-AC-OPF problems. The network sizes for these systems ranged from 500 to 30,000 buses (transmission network nodes)~\citep{GoCompetition}. The operating cost and constraint violation penalties achieved by our algorithm resulted in a second-place finish overall among the 27 teams in the GO Challenge 1 competition. Our algorithm obtained the best scores for 20\% of the problems, the best or second-best scores for 53\% of the problems, and top-ten scores for almost all ($>99\%$) of the problems. The results of this competition, by our team and others, indicate that nonlinear optimization techniques such as those used in the algorithm described in this paper are effective at solving industrially relevant, large-scale SC-AC-OPF problems.

This paper is organized as follows. Section~\ref{sec.problem} introduces the problem formulation and discusses several key modeling decisions. Section~\ref{sec.algorithm} describes our proposed solution algorithm, including an overall decomposition strategy along with our contingency selection and evaluation procedures.  Section~\ref{sec.results} provides the results of numerical experiments that isolate the effectiveness of various components of our overall algorithm. Section~\ref{sec.conclusion} provides concluding remarks.

\section{Problem Statement}\label{sec.problem}

A complete description of the problem statement can be found on the GO Competition website \citep{GoCompetition}.  In this section, we present a somewhat simplified statement of the problem along with a high-level model through which we can describe our decomposition strategy.  We also discuss the most salient modeling decisions we made that we believe benefit our algorithm's performance.

\bfigure
  \centering
  \includegraphics[width=0.8\textwidth,clip=true,trim=20 60 80 120]{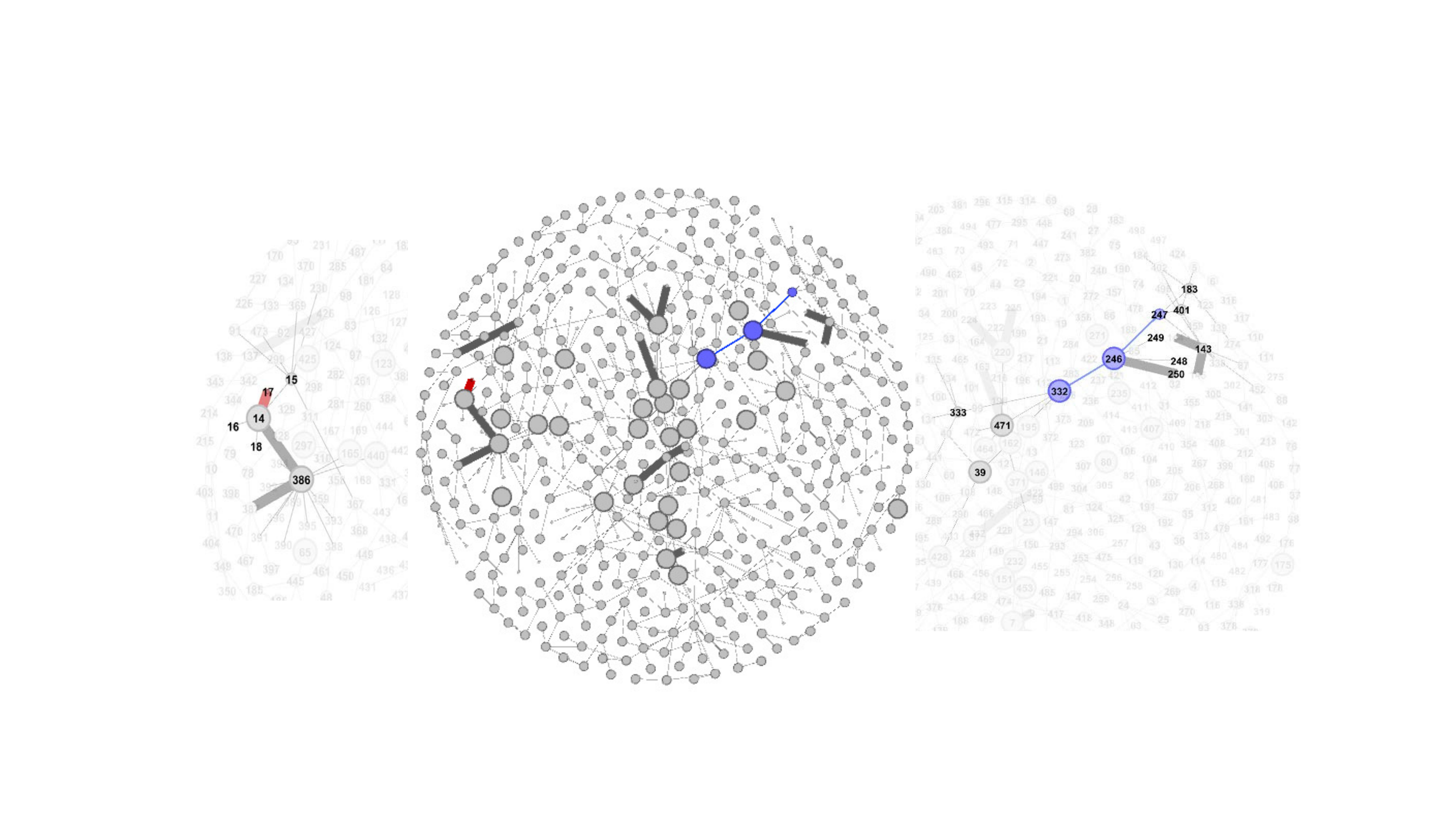}
  \caption{Visualization of Network 1 with 500 nodes. The thick lines represent parallel lines or transformers and the size of each node indicates the voltage of the bus. Bus 17 (see the red edge, with a close-up view shown on the left) is associated with the generator contingency that has the highest penalty value. Transformer 247 to 246 and line 246 to 332 (both in blue, with a close-up view shown on the right) are the transformer and line contingencies associated with the highest penalty values. All three of the indicated contingencies are located near a parallel line and close to a high-voltage bus.}
\efigure

\subsection{Power flow model}\label{sec.base_case}

Generally speaking, we follow the notation from the official problem statement of the GO Challenge~1 competition, with certain simplifications to ignore details that are not of consequence for the main ideas of this paper.  As is typical in power system analyses, we adopt the \emph{per-unit} system, with the input data provided by the competition normalized accordingly.  As a result, instead of actual voltage, current, and power units (volts, amps, and watts), we represent these quantities in our model by normalized per-unit quantities.

Let $\Ical$ be the set of buses, $\Gcal$ be the set of generators, $\Ecal$ be the set of lines, and $\Fcal$ be the set of transformers.  (For our purposes throughout the paper, we ignore the distinction of \emph{areas} defined in the official problem statement.)  Bus $i \in \Ical$ has as decision variables its voltage magnitude $v_i$, voltage angle $\theta_i$, and controllable shunt susceptance $b_i^{CS}$.  Generator $g \in \Gcal$ has as decision variables its active power output $p_g$ and reactive power output~$q_g$.

The flow on each line is defined by AC power flow for a $\Pi$-circuit model.  Specifically, corresponding to line $e \in \Ecal$, the active and reactive power flows into the line from an \emph{origin} bus $i_e^o \in \Ical$ to a \emph{destination} bus $i_e^d \in \Ical$ are, respectively,
\bequation\label{eq.line_flow}
  \baligned
    p_e^o
      =&\ g_e v_{i_e^o}^2 - (g_e \cos(\theta_{i_e^o} - \theta_{i_e^d}) + b_e \sin(\theta_{i_e^o} - \theta_{i_e^d})) v_{i_e^o} v_{i_e^d} \\
    q_e^o
      =&\ -(b_e + b_e^{CH}/2) v_{i_e^o}^2 + (b_e \cos(\theta_{i_e^o} - \theta_{i_e^d}) - g_e \sin(\theta_{i_e^o} - \theta_{i_e^d})) v_{i_e^o} v_{i_e^d},
  \ealigned
\eequation
where $g_e$ is the conductance, $b_e$ is the susceptance, and $b_e^{CH}$ is the total charging susceptance for the line.  The active and reactive power flows into the line from the destination bus to the origin bus are defined similarly, with the roles of the origin and destination buses reversed.  Also similarly, corresponding to a transformer $f \in \Fcal$, the active and reactive power flows into the transformer from an origin bus $i_f^o \in \Ical$ to a destination bus $i_f^d \in \Ical$ are defined respectively by
\bequation\label{eq.transformer_flow}
  \baligned
    p_f^o
      =&\ (g_f/\tau_f^2 + g_f^M) v_{i_f^o}^2 -(g_f/\tau_f \cos(\theta_{i_f^o} - \theta_{i_f^d} - \theta_f) + b_f/\tau_f \sin(\theta_{i_f^o} - \theta_{i_f^d} - \theta_f)) v_{i_f^o} v_{i_f^d} \\
    q_f^o
      =&\ -(b_f/\tau_f^2 + b_f^M) v_{i_f^o}^2 + (b_f/\tau_f \cos(\theta_{i_f^o} - \theta_{i_f^d} - \theta_f) - g_f/\tau_f \sin(\theta_{i_f^o} - \theta_{i_f^d} - \theta_f)) v_{i_f^o} v_{i_f^d},
  \ealigned
\eequation
where $g_f$ is the conductance, $\tau_f$ is the tap ratio, $\theta_f$ is the phase shift, $g_f^M$ is the magnetizing conductance, $b_f$ is the susceptance, and $b_f^M$ is the magnetizing susceptance for the transformer.  The active and reactive power flows into the transformer from the destination bus are defined similarly, with the origin/destination roles reversed and slightly different coefficients; for details, see the official problem formulation.  Note that all parameters in both the line and transformer models are specified constants, including the transformer phase shifts $\{\theta_f\}_{f\in\Fcal}$ and tap ratios $\{\tau_f\}_{f\in\Fcal}$.

Power balance equations ensure that, at each bus $i \in \Ical$, the active and reactive power outputs from generators are equal to the sums of flows into other components at the bus.  In the model, such balance equations are relaxed with nonnegative slack variables that are penalized in the objective function.  The resulting relaxed balance equations are, respectively,
\bequation\label{eq.bus_balance}
  \baligned
    \sum_{g \in \Gcal_i} p_g - p_i^L - g_i^{FS}v_i^2 - \sum_{e \in \Ecal_i^o} p_e^o - \sum_{e \in \Ecal_i^d} p_e^d - \sum_{f \in \Fcal_i^o} p_f^o - \sum_{f \in \Fcal_i^d} p_f^d &= \sigma_i^{P_+} - \sigma_i^{P_-} \\
    \sum_{g \in \Gcal_i} q_g - q_i^L + (b_i^{FS} + b_i^{CS})v_i^2 - \sum_{e \in \Ecal_i^o} q_e^o - \sum_{e \in \Ecal_i^d} q_e^d - \sum_{f \in \Fcal_i^o} q_f^o - \sum_{f \in \Fcal_i^d} q_f^d &= \sigma_i^{Q_+} - \sigma_i^{Q_-}.
  \ealigned
\eequation
Here, for bus $i \in \Ical$, $\Gcal_i \subseteq \Gcal$ is the subset of generators connected to the bus, $\Ecal_i^o \subseteq \Ecal$ (resp.,~$\Ecal_i^d  \subseteq \Ecal$) is the subset of lines with origin (resp.,~destination) as the bus, $\Fcal_i^o \subseteq \Fcal$ (resp.,~$\Fcal_i^d \subseteq \Fcal$) is the subset of transformers with origin (resp.,~destination) as the bus, $p_i^L$ (resp.,~$q_i^L$) is the bus' constant active (resp.,~reactive) power load, $g_i^{FS}$ (resp.,~$b_i^{FS}$) is the bus' fixed shunt conductance (resp.,~susceptance), and $\{\sigma_i^{P_+},\sigma_i^{P_-},\sigma_i^{Q_+},\sigma_i^{Q_-}\}$ is a set of nonnegative slack variables.

For each line $e \in \Ecal$, the current flow is penalized beyond a bound defined by the voltage and a current limit $\overline{R}_e$.  Specifically, at the origin and destination of the line, the constraints are given, with a common nonnegative slack variable $\sigma_e^S$, respectively by
\bequation\label{eq.line_current_rating}
  \baligned
    (p_e^o)^2 + (q_e^o)^2 &\leq (\overline{R}_e v_{i_e^o} + \sigma_e^S)^2 \\
    (p_e^d)^2 + (q_e^d)^2 &\leq (\overline{R}_e v_{i_e^d} + \sigma_e^S)^2. \\
  \ealigned
\eequation
For each transformer, the power flow is penalized in a similar manner beyond an upper bound defined by the transformer's apparent power limit $\overline{s}_f$.  Specifically, at the origin and destination of the line, the constraints are given, with a common nonnegative slack variable $\sigma_f^S$, respectively by
\bequation\label{eq.transformer_power_rating}
  \baligned
    (p_f^o)^2 + (q_f^o)^2 &\leq (\overline{s}_f + \sigma_f^S)^2 \\
    (p_e^d)^2 + (q_e^d)^2 &\leq (\overline{s}_f + \sigma_f^S)^2. \\
  \ealigned
\eequation

\subsection{Contingency response}\label{sec.cont_case}

Each problem instance comes with a set of contingencies, indexed over a set~$\Kcal$, the models for which link to the \emph{base case} variables and constraints defined in~Section~\ref{sec.base_case}.  All variables and constraints are duplicated for each contingency, except for a missing generator, line, or transformer.  To respond to the loss of a component in each contingency, active and reactive power generation at each available generator can be adjusted from the base case value according to specified policies.  For the sake of brevity, we do not repeat all of the constraints pertaining to the power flow model for each contingency.  Instead, we present only the constraints that define the contingency response.

Let $\Gcal_k \subseteq \Gcal$ be the subset of generators that are available in contingency $k \in \Kcal$, and let $\Gcal_k^P \subseteq \Gcal_k$ be the subset of those generators that are chosen to respond in the contingency.  For generators that are available, but are chosen not to respond, the active power output in the contingency~$p_{gk}$ equals the base case value~$p_g$.  Otherwise, the active power output is adjusted by a common perturbation~$\Delta_k$, leading to a complementarity constraint that at least one of two sets of conditions holds (where $\underline{p}_g$ and $\overline{p}_g$ are, respectively, lower and upper bounds on the active power output):
\bequation\label{eq.comp_active}
  \baligned
    0 \leq (p_{gk} - \underline{p}_g) &\perp (p_g + \alpha_g \Delta_k - p_{gk}) \leq 0 \\ \text{or}\ \ 
    0 \geq (p_{gk} - \overline{p}_g) &\perp (p_g + \alpha_g \Delta_k - p_{gk}) \geq 0.
  \ealigned
\eequation
(Here, for scalars $\eta$ and $\mu$, the expression $\eta \perp \mu$ means that either $\eta$ or $\mu$ is equal to zero.)  Similarly for reactive power, with the reactive power output in the contingency denoted as $q_{gk}$, the contingency response constraints are given for all $g \in \Gcal_k$ (where $\underline{q}_g$ and $\overline{q}_g$ are, respectively, lower and upper bounds on the reactive power output) by
\bequation\label{eq.comp_reactive}
  \baligned
    0 \leq (q_{gk} - \underline{q}_g) &\perp (v_{i_g} - v_{i_gk}) \leq 0 \\ \text{or}\ \ 
    0 \geq (q_{gk} - \overline{q}_g) &\perp (v_{i_g} - v_{i_gk}) \geq 0.
  \ealigned
\eequation
These constraints model a \emph{preventative} formulation, as opposed to a \emph{corrective} formulation that would allow for more general responses to each contingency.  The active power response models the effect of commonly used frequency control schemes~\citep{jaleeli1992}, while the reactive power response models so-called ``PV/PQ switching''~\citep{stott1974}.

\subsection{High-level model}

At a high level, the problem can be expressed in the manner shown in \eqref{prob.high_level}.  As defined in more detail below, let $u_0$ be the control variables in the base case, $y_0$ be the state variables in the base case, and $(\varsigma_0,\varsigma_0^+,\varsigma_0^-)$ be the slack variables in the base case.  Similarly, for each contingency $k \in \Kcal$, let $u_k$ be the control variables, $y_k$ be the state variables, and $(\varsigma_k,\varsigma_k^+,\varsigma_k^-)$ be the slack variables.

\bfigure[ht]
\bequation\label{prob.high_level}
  \baligned
    \min_{\substack{(u_0,y_0,\varsigma_0,\varsigma_0^+,\varsigma_0^-) \\ \{(u_k,y_k,\varsigma_k,\varsigma_k^+,\varsigma_k^-)\}_{k\in K}}} &\ f_0(u_0,y_0) + \psi_0(\varsigma_0,\varsigma_0^+,\varsigma_0^-) + \frac{1}{|\Kcal|} \sum_{k\in\Kcal} \psi_k(\varsigma_k,\varsigma_k^+,\varsigma_k^-) \\
    \st &\ \left\{
    \baligned
      a_0(u_0,y_0) &= 0 \\
      b_0(u_0,y_0) &= \varsigma_0^+ - \varsigma_0^- \\
      c_0(u_0,y_0) &\leq \varsigma_0 \\
      \underline{u}_0 &\leq u_0 \leq \overline{u}_0 \\
      \underline{y}_0 &\leq y_0 \leq \overline{y}_0 \\
      (\varsigma_0,\varsigma_0^+,\varsigma_0^-) &\geq 0 \\
      a_k(u_k,y_k) &= 0 && \text{for all}\ k \in K \\
      b_k(u_k,y_k) &= \varsigma_k^+ - \varsigma_k^- && \text{for all}\ k \in K \\
      c_k(u_k,y_k) &\leq \varsigma_k && \text{for all}\ k \in K \\
      \underline{u}_k &\leq u_k \leq \overline{u}_k && \text{for all}\ k \in K \\
      \underline{y}_k &\leq y_k \leq \overline{y}_k && \text{for all}\ k \in K \\
      (\varsigma_k,\varsigma_k^+,\varsigma_k^-) &\geq 0 && \text{for all}\ k \in K \\
      0 \leq \tau_k^L(u_0,y_0,u_k,y_k) \perp \tau_k^R(u_k,y_k) &\geq 0 && \text{for all}\ k \in K.
    \ealigned
    \right.
  \ealigned
\eequation
\efigure

At a more detailed level, the base case control variables are the active and reactive power produced by the generators and the controllable shunt susceptances at the buses, whereas the only control variables in the contingencies are the controllable shunt susceptances at the buses.  All other variables are state variables except those that are specified as slack variables that are penalized in the objective function.  The objective consists of a base case objective $f_0$---representing the cost of active power generation---along with penalty terms for the base case ($\psi_0$) and contingency slack variables ($\{\psi_k\}_{k\in\Kcal}$).  The base case objective is a convex piecewise linear cost of generation, the definition of which involves some auxiliary variables; see the official problem statement.

The constraints come in a few different types.  Besides simple bounds on the control, state, and penalty variables, they are as follows.  The constraint functions $a_0$ and $\{a_k\}_{k\in\Kcal}$ capture the line and transformer flow definitions \eqref{eq.line_flow}--\eqref{eq.transformer_flow}, where $a_0$ fixes a single phase angle in the base case (to zero) in order to set a reference angle.  The constraint functions $b_0$ and $\{b_k\}_{k \in \Kcal}$ capture the bus power balance \eqref{eq.bus_balance}.  The constraint functions $c_0$ and $\{c_k\}_{k \in \Kcal}$ capture the line current and transformer power ratings constraints \eqref{eq.line_current_rating}--\eqref{eq.transformer_power_rating}.  Finally, the constraint functions $\{\tau_k^L\}_{k \in \Kcal}$ and $\{\tau_k^R\}_{k \in \Kcal}$ capture the complementarity constraints \eqref{eq.comp_active}--\eqref{eq.comp_reactive}.

It is worth emphasizing that only the complementarity constraints involve a combination of base case and contingency variables, and \emph{no} single constraint involves variables from more than one contingency at a time.  This makes the problem nearly separable, which is a property that needs to be exploited in a solution approach for it to be efficient.

\subsection{Modeling decisions}

The aforementioned formulation allows some flexibility, and the specific modeling choices that are made can have a significant impact on the performance of a solution algorithm.  In this subsection, we comment on the modeling choices that we believe are most consequential in our approach.

First, as is common when solving many types of optimization problems, one has the option of eliminating a potentially large number of constraints.  For example, in the problem at hand, one could eliminate the line and transformer flow definition constraints \eqref{eq.line_flow}--\eqref{eq.transformer_flow} and simply plug the expressions for the variables being defined by these constraints into the other constraints.  This would have the effect of eliminating a number of variables and number of constraints on the order of the number of lines plus the number of transformers, multiplied by the number of contingencies.  However, the downside of such elimination is a large increase in the density of the constraint Jacobians.  Hence, rather than perform this elimination explicitly, we formulate the problem with all of these variables and constraints and allow the linear system solver in the nonlinear optimization method to exploit the sparsity of the resulting linear systems.

Second, the line current and transformer power ratings constraints \eqref{eq.line_current_rating}--\eqref{eq.transformer_power_rating} are essentially upper bounds on the norms of two-dimensional vectors involving active and reactive power at an origin or destination bus.  The norm is a nonsmooth function, which might be problematic for an interior-point method whose theoretical guarantees depend on smoothness of the problem functions.  Hence, the constraints that we state in \eqref{eq.line_current_rating}--\eqref{eq.transformer_power_rating} are squared versions of the constraints stated in the official problem statement so that $c_0$ and $\{c_k\}_{k\in\Kcal}$ are smooth.

Third, the generator active and reactive power contingency responses can be modeled in various ways.  For example, the official problem formulation document provided by \mbox{ARPA-E} includes three possible formulations: a \emph{logical formulation} in which the responses are captured by complementarity constraints, a \emph{projection formulation}, and a \emph{mixed integer programming} formulation.  We saw little benefit in the projection formulation over the logical formulation, and saw little hope in being able to provide good solutions in a limited timeframe using mixed-integer programming techniques.  Hence, as previously indicated, we choose to follow the logical formulation.

Finally, we note that the power flow equations~\eqref{eq.line_flow} and~\eqref{eq.transformer_flow} are based on voltage phasors in polar coordinates, i.e., voltage magnitudes $v$ and angles $\theta$. These equations could be rewritten in a variety of ways~\citep{ferc1,molzahn_hiskens-fnt2019}. For instance, the power flow equations~\eqref{eq.line_flow} and~\eqref{eq.transformer_flow} can be represented as systems of quadratic polynomials through the use of rectangular voltage coordinates $v_{d,i} = v_i\cos(\theta_i)$ and $v_{q,i} = v_i\sin(\theta_i)$.  After some initial experiments through which we considered alternative formulations, we ended up following the polar formulation from the \mbox{ARPA-E} problem description.  Building on studies such as that in~\citep{ferc5}, it would be worthwhile to evaluate these alternatives more comprehensively.

\subsection{Competition Requirements}

An implementation of an algorithm submitted for the GO Challenge 1 competition to solve problems of the form in~(\ref{prob.high_level}) had to adhere to strict guidelines dictated by the competition organizers.  These guidelines strongly influenced our algorithmic development.  Competitors were evaluated within four divisions, which were defined by different scoring methodologies and strict time limits.

The year-long competition involved three trial events before the final submission deadline, with increasingly difficult datasets released before and after each trial event for testing and evaluation.  These datasets consisted of input files describing the network and scenarios provided by the competition, formatted according to power system industry standards.  For evaluation and scoring purposes, the competition had two scoring methods, one focusing on the objective value offered by the final solution and another focusing on the robustness of the algorithm relative to the other teams' submissions.  Both the generation costs associated with power production and any penalty terms (on the slack variables) were reflected in the final objective value, thus ensuring feasibility of each team's outputs, despite any potential constraint violations.

Algorithms were required to provide a solution for the base case variables defined in Section~\ref{sec.base_case} as well as solutions for each of the contingency scenario variables defined in Section~\ref{sec.cont_case}. The submitted software was divided into two components: \emph{Code 1}, which determined the base case solution, i.e., values of the variables associated with nominal operation, and \emph{Code 2}, which produced the contingency response solutions for every contingency scenario.  More precisely, Code 1 observed the combined base case and contingencies problem in order to optimize the base case solution subject to the contingencies, but was only asked to provide a base case solution.  Afterwards, Code 2 was run to produce contingency response solutions using the base case solution that was determined by the run of Code 1.  The competition was also divided into two categories: real-time and offline optimization with strict 10-minute and 45-minute time limits, respectively, imposed on the execution of \emph{Code 1}.  An additional time limit of two seconds per contingency was imposed on \emph{Code 2} regardless of category.

The efficiency of the algorithm and the speed of the software were very important to ensure that our computations were performed within the prescribed time limits.  In addition to time limits, the scoring methodologies and solution output requirements also informed many of our modeling, algorithmic, and implementation choices.  Chief among these was the computation of the base case variables in the first stage by iteratively identifying the most \emph{important} contingencies and only explicitly imposing the corresponding constraints.  (We explain further our meaning of the \emph{importance} of a contingency in Section~\ref{sec.algorithm}.)  This approach traded off (ideally small) penalties associated with the unconsidered contingencies in return for computational speed improvements. Additionally, both the first and second stages of our algorithm began by applying a fast, but approximate contingency evaluation procedure in order to avoid the computationally intensive task of fully evaluating the penalties associated with each contingency.  Moreover, we imposed computational limitations on many of the steps in our algorithm in order to make the best use of the total available time.  We explain further these features of the competition and their implications for our algorithmic developments in Section~\ref{sec.algorithm}.

\section{Algorithm}\label{sec.algorithm}

A high-level view of our algorithm is presented in Figure~\ref{fig.flow_chart}.  After a preprocessing phase (about which we leave discussion until Section~\ref{sec.preprocessing}), a base case solution is obtained by solving \eqref{prob.high_level} with \emph{no} contingencies.  Given this base case solution, a scheme is employed to determine an initial ranking of contingencies according to their estimated importance (see Section~\ref{sec.initial_ranking}), after which a parallel process is initiated to start providing better estimates of each contingency's importance (see Section~\ref{sec.contingency_evaluation}).  This process involves a specialized procedure for handling the complementarity constraints (see Section~\ref{sec.complementarities}) and one for fast (approximate) contingency evaluation (see Section~\ref{sec.fast}).  At this point, a few contingencies are selected and the algorithm enters an iterative process during which a \emph{master problem} is solved---i.e., a problem of the same form as \eqref{prob.high_level}, but with only those contingencies that have been selected as important---and contingencies are continually evaluated and possibly selected for inclusion in the master problem in subsequent iterations.  Many of the procedures in this loop are performed in a coordinated fashion in parallel.  We also mention in this section the techniques that we employed for avoiding the selection of \emph{dominated} contingencies (see Section~\ref{sec.dominated}) as well as a summary of our entire algorithm (see  Section~\ref{sec.overall}).

\bfigure[ht]
  \centering
  \tikzstyle{block} = [rectangle, draw, text width=12em, text centered, rounded corners, minimum height=4em, font=\large]
  \tikzstyle{par_block} = [rectangle, draw, fill=blue!20, text width=12em, text centered, rounded corners, minimum height=4em, font=\large]
  \tikzstyle{line} = [draw, -latex']
  \tikzstyle{cloud} = [draw, ellipse,fill=red!20, node distance=3cm, minimum height=2em]
  \resizebox{.8\textwidth}{!}{%
  \begin{tikzpicture}[node distance = 2cm, auto, trans/.style={,<->,shorten >=2pt,shorten <=2pt,>=stealth}]
    \node [block] (preprocessing) {Preprocessing};
    \node [block, below of=preprocessing] (base) {Base case\\ solution: Ipopt}; 
    \node [block, below of=base] (write1) {Write Solution File};
    \node [block, right of=write1, xshift=22cm] (pre2) {Preprocessing}; 
    \node [block, below of=write1] (ranking) {\mbox{Initial Contingency} Ranking};
    \node [block, below of=ranking] (selection_init) {Contingency Selection};
    \node [par_block, below of=selection_init] (master) {Master\\ Problem Solve};
    \node [block, below of=master] (write2) {Write Solution File};
    \node [par_block, below of=write2] (selection){ Contingency Selection};
    \node[par_block, right of= ranking, xshift=10cm, yshift=-2cm](scheduler){\mbox{Parallel Parent} Scheduler};
    \node [par_block, below of=scheduler, xshift=-4cm, yshift=-3cm] (evaluation_1){Contingency Evaluation};
    \node[block, below of=pre2](read){Read Solution File};
    \node[par_block,below of=read](evaluation){Contingency Evaluation};
    \node [par_block, below of=scheduler, xshift=4cm,yshift=-3cm] (evaluation_3){ Contingency Evaluation};
    \node[] at (-4.2,-5) {\Large \emph{Code 1}};
    \node[] at (22,-14.5) {\Large \emph{Code 2}};
    \background{evaluation_1}{pre2}{evaluation}{evaluation_1}{code2}
    \background{preprocessing}{preprocessing}{evaluation_3}{selection}{code1}
    \path [line] (preprocessing) -- (base);
    \path [line] (base) --(write1);
    \path[line] (write1)--(ranking);
    \path [line] (ranking) -- (selection_init);
    \path [line] (selection_init) -- (master);
    \draw [line] (master) -- (write2);
    \draw[trans] (write2)--(selection);
    \draw[line] (pre2)--(read);
    \draw[line] (read)--(evaluation);
    \draw [thick,loosely dotted] (evaluation_1)--(evaluation_3);
    \draw[trans](scheduler) --(evaluation_1);
    \draw[trans](scheduler) --(evaluation_3);
    \draw[trans](selection_init) -- (scheduler);
    \draw[trans](selection) -- (scheduler);
    \draw[trans](evaluation)--(scheduler);
  \end{tikzpicture}
  } 
  \caption{High-level view of algorithm.}
  \label{fig.flow_chart}
\efigure

Overall, our algorithmic strategy is to identify quickly what are the \emph{important} contingencies to include in the master problem in order to obtain a solution that is as close as possible to the solution that would be obtained if all contingencies were to be included.  By the importance of a contingency, we mean its effect on the optimal value of the overall problem.  Generally speaking, if an important contingency is not included in the master problem, then this means that the solution obtained for the master problem would represent a poor solution for the overall problem, whereas if a less important contingency is ignored, then the solution would not be much poorer.  A complicating factor is that, in some cases, two contingencies could be individually important, although as soon as one is included in the master problem, the other becomes less important.  For example, this may be the case for contingencies corresponding to parallel lines in the transmission network.  Our overall strategy attempts to determine as few contingencies as possible such that, once these are included in the master problem, the remaining contingencies are unimportant. We note that similar algorithmic strategies are used in previous power systems literature; see, e.g., \citep{wollenberg1979,stott1987,capitanescu2016}.

All optimization problems in the algorithm---including the instances of the master problem and penalty minimization problems for contingency evaluation---are solved with the interior-point method implemented in the Ipopt software~\citep{WaecBieg06}. Interior-point methods were first employed to solve power system problems in the early 1990s for the purpose of state estimation~\citep{clements1991}. Subsequent work over the next three decades applied interior-point methods to solve various OPF problems; see, e.g., \citep{wu1994,wang2007,capitanescu2011,capitanescu2013}.


\subsection{Initial Ranking of Contingencies}\label{sec.initial_ranking}

For a given base case solution, evaluation of the penalty associated with contingency $k \in \Kcal$ is relatively expensive computationally, since it involves solving a nonlinear optimization problem that contains all of contingency $k$'s variables; see Section~\ref{sec.contingency_evaluation}.  Hence, after a base case solution is obtained, it may be detrimental to wait until all contingencies have been evaluated in this manner before solving a new master problem involving a few contingencies to obtain an improved base case solution.  For this reason, we employ a contingency ranking heuristic in order to more quickly identify potentially important contingencies. These heuristics are deployed immediately after we obtain the first base case solution. The contingencies that are identified as potentially being the most important are evaluated first using the contingency evaluation and selection strategies that we describe in subsequent subsections. Through this process, which exploits parallel computations, not all contingencies are evaluated before a new master problem is solved.



For several decades, contingency ranking has been an extensively studied subject in the power systems literature (see, e.g,.~\citep{wollenberg1979} and~\citep{stott1987} as well as the more recent survey in~\citep{wu2017}). While inspired by this literature, we created a new contingency ranking heuristic based on the problem formulation, dataset characteristics, and computing time requirements of the GO competition.

Our contingency ranking heuristic uses features of the transmission network topology as well as the values obtained from the initial base case solution. To develop this heuristic, we first identified a collection of potential features, such as topology information in the neighborhood of a contingency and generation or transmission flow loss introduced by the contingency. We then trained a ridge regression 
model, i.e., a least squares loss function with an $\ell_2$-norm regularization term, using a linear combination of these features to predict the penalty associated with each contingency---whether it be defined by a generator, line, or transformer failure---in any scenario using the datasets provided by the competition. The dataset was divided into training and testing groups. The model parameters were trained on the training dataset and tested on the rest of the networks including both new scenarios for the networks in the training dataset as well as new networks that were not included in the training set. We repeated the cross validation process 10 times and took averages across these repetitions to get the final model parameters. During the execution of our algorithm, we use the obtained model to generate penalty values that predict the importance of each contingency. Sorting the contingencies in descending order of their penalty values provides our initial contingency ranking. We next describe the features we use to generate this ranking.

Studying the testing datasets offered by the competition organizers suggests that there is often a direct correlation between the importance of a contingency and the degree (i.e., number of neighboring buses) of the buses near the contingency. We consider the voltage ratings of either the bus(es) associated with a contingency (i.e., the bus where the generator contingency occurs or the end points of either line or transformer contingencies) or the \emph{neighboring} buses, with the idea that a generator, line, or transformer contingency may be more important if one of these ratings is sufficiently large. The power flow losses introduced by the contingencies are also identified as important. Specifically, after considering many combinations of potential features, for each contingency $k \in \Kcal$, the best predicted importance is given by a linear combination of the following features: $\{ t^g_k, t^l_k, t^t_k, l^p_k, l^s_k, l^c_k, v^d_k, d^o_k, d^d_k, \pi_k\}$.
The first three features, namely, $\{t^g_k, t^t_k, t^l_k\}$ are binary variables indicating the type of contingency: for a generator contingency, $t^g_k=1$ and $t^t_k=t^l_k=0$; for a line contingency, $t^l_k=1$ and $t^g_k=t^t_k=0$; and for a transformer contingency, $t^t_k=1$ and $t^g_k=t^l_k=0$.  The remainder of the features have different interpretations  for generator versus line or transformer contingencies, as described next.

For contingency $k \in \Kcal$ associated with generator $g$, we define three features related to the loss in power generation induced by the generator failure relative to the base case solution: active power generation loss  $l^p_k = p_g$, apparent power generation loss $l^s_k = \sqrt{p_g^2+q_g^2}$, and apparent power loss relative to the generator capacity $l^c_k = \frac{\sqrt{p_g^2+q_g^2}}{\sqrt{\overline{p}_g^2+\overline{q}_g^2}}$. We also define the following two features related to the neighboring buses of the generator: $v^d_k$ is the highest voltage rating prior to per-unit normalization (given by the field BASKV in the dataset) of any neighboring buses and $d^o_k$ is the degree of the bus where the generator $g$ is located. The rest of the features are set to zero (i.e., $d^t_k = \pi_k = 0$) for all generator contingencies. 

For contingency $k \in \Kcal$ associated with line $e$ or transformer $f$ connecting origin bus $o$ and destination bus $d$, we similarly define three features related to the flow loss relative to the base case solution.  In particular, we define the active power flow loss as $l^p_k = \max\{p_e^o, p_e^d\}$ for line $e$ (and similarly for transformer $f$), the apparent power loss as $l^s_k =  \max\big\{\sqrt{(p_e^o)^2 + (q_e^o)^2 }, \sqrt{ (p_e^d)^2 + (q_e^d)^2}\big\}$ for line $e$ (and similarly for transformer $f$), and the apparent power loss relative to the flow capacity as $l^c_k = \max\left\{\sqrt{\frac{(p_e^o)^2 + (q_e^o)^2}{\overline{R}_e v_{i_e^o}}}, \sqrt{\frac{(p_e^d)^2 + (q_e^d)^2}{\overline{R}_e v_{i_e^d}}}\right\}$ for line $e$ and  $l^c_k = \max\left\{\sqrt{\frac{(p_e^o)^2 + (q_e^o)^2}{\bar s_f}}, \sqrt{\frac{(p_e^d)^2 + (q_e^d)^2}{\bar s_f}}\right\}$ for transformer $f$. The following four features related to the network topology also enhance the quality of predicted penalty values: $v_k^d$, the voltage rating of bus $d$; $d_k^o$, the degree of bus $o$; $d_k^d$, the degree of bus $d$; and $\pi_k$, a weight set to $10$ if the line/transformer has a parallel counterpart, i.e., there exists another line/transformer connecting the same pair of origin and destination buses, or set to $0$ otherwise. The model's performance does not improve when we include additional features such as the degree of the origin and destination buses for the lines/transformers.

Among these features, the apparent power loss $l^s_k$ has the largest coefficient. This observation agrees with our numerical studies in Section \ref{sec.num_initial_ranking}, where we show that this feature alone has very good predictive power.  One interesting property of our model is that it is universal to all networks regardless of their sizes. Introducing network size as a feature does not improve the predictive power of the model. We suspect this is due to the effect of a contingency often being either local or well modeled by the power loss features alone.

We remark in passing that the aforementioned model uses a simple linear combination of the features. Alternative approaches (e.g., multiplicatively combining the features, using additional features, etc.) generally led to similar or inferior empirical performance in our tests. Finally, we note that some contingencies are often identified as important (in terms of their typical penalty values) across many scenarios for the same system. For each scenario for a particular system, we recorded the top contingencies which incurred a large penalty value. We then took the union of all of these top contingencies to form a candidate list.  The initial contingency ranking in a run of our algorithm combines this candidate list with the generator, line, and transformer ranking values described above to create an initial prioritized list of the contingencies.

\subsection{Contingency Evaluation}\label{sec.contingency_evaluation}

Given a base case solution $(u_0,y_0)$, a solution for contingency $k \in \Kcal$ is obtained by solving the (continuous and nonlinear) penalty minimization problem

\bequation\label{eq:slow_evaluation}
  \baligned
    \min_{u_k,y_k,\varsigma_k,\varsigma_k^+,\varsigma_k^-} &\ \psi_k(\varsigma_k,\varsigma_k^+,\varsigma_k^-) \\
    \st &\ \left\{
    \baligned
      a_k(u_k,y_k) &= 0 \\
      b_k(u_k,y_k) &= \varsigma_k^+ - \varsigma_k^- \\
      c_k(u_k,y_k) &\leq \varsigma_k \\
      \underline{u}_k \leq u_k &\leq \overline{u}_k \\
      \underline{y}_k \leq y_k &\leq \overline{y}_k \\
      (\varsigma_k,\varsigma_k^+,\varsigma_k^-) &\geq 0 \\
      0 \leq \tau_k^L(u_0,y_0,u_k,y_k) \perp \tau_k^R(u_k,y_k) &\geq 0,
    \ealigned
    \right.
  \ealigned
\eequation
where the variables and problem functions are defined as in \eqref{prob.high_level}.  We note that one might consider warm-starting the solve for \eqref{eq:slow_evaluation}, e.g., using the base case solution values or a solution of the problem obtained in a previous iteration of the overall algorithm.  However, we did not find warm-starting to work well in our experiments, which may be due to the well-known difficulty inherent in effectively warm-starting interior-point methods.

The complementarity constraints in \eqref{eq:slow_evaluation} are not included as explicit complementarity constraints when the problem is solved in our algorithm.  Rather, they are replaced by a set of (smooth) constraints that depend on an active-set prediction regarding which of the terms involved in the complementarity constraints are active and which are free.  Overall, in our solution algorithm, \eqref{eq:slow_evaluation} is solved approximately through a loop in which the complementarity predictions are updated iteratively.  We describe this process in further detail in the next subsection.

\subsection{Handling Complementarities}\label{sec.complementarities}

Each pair of complementarity constraints in \eqref{eq.comp_active} and \eqref{eq.comp_reactive} involves a variable, call it $\chi$ (representing $p_{gk}$ in \eqref{eq.comp_active} and $q_{gk}$ in \eqref{eq.comp_reactive}), and a linear expression, call it $\rho$ (representing $p_g + \alpha_g \Delta_k - p_{gk}$ in \eqref{eq.comp_active} and $v_{i_g} - v_{i_gk}$ in \eqref{eq.comp_reactive}).  Each constraint requires that at least one of the following holds:
\bitemize[label={$\bullet$}]
  \item $\chi$ equals a lower bound $\underline\chi$ and $\rho \leq 0$,
  \item $\chi$ equals an upper bound $\overline\chi$ and $\rho \geq 0$, or
  \item $\chi \in [\underline\chi,\overline\chi]$ and $\rho = 0$.
\eitemize
Our strategy for handling complementarities may be referred to as an active set approach, wherein we iteratively (i) make a prediction about which of these conditions holds at an optimal solution, (ii) solve the resulting problem (that involves no complementarity constraints), then (iii) update our prediction based on multiplier values obtained from the solution obtained in step (ii).

Let us refer to the condition $\chi = \underline\chi$ as the \emph{lower} segment, $\chi \in [\underline\chi,\overline\chi]$ as the \emph{middle} segment, and $\chi = \overline\chi$ as the \emph{upper} segment of a complementarity constraint.  The corresponding restrictions on the linear expression are $\rho \leq 0$, $\rho = 0$, and $\rho \geq 0$, respectively.  Our procedure for updating a prediction for a complementarity is based on the same strategy in all cases.  For example, let us consider the case when the prediction in step (i) is that the lower segment is optimal.  In step (ii), the problem (without complementarities) is solved with the constraints $\chi = \underline\chi$, $\rho - s = 0$, and $s \leq 0$, where $s$ is a slack variable.  Let $\lambda \geq 0$ be the Lagrange multiplier for $s \leq 0$.  If $\lambda > 0$ and $-s/\lambda < 10^{-6}$, then the constraint $s \leq 0$ is deemed to be active and the prediction is altered to say that the middle segment is optimal; otherwise, the prediction is not changed.  Using such a procedure, a prediction may also be altered from middle to lower, middle to upper, or upper to middle.

Our active-set estimates for the complementarity constraints are updated in this manner during contingency evaluation only.  In other words, the active-set estimates are not updated during a master problem solve.  We considered updating the predictions iteratively for the master problem as well, but this turned out to be too time-consuming due to the greater expense of solving each master problem, and might not be worthwhile anyway since the predictions are updated during the next contingency evaluations in any case.  During contingency evaluation, the predictions are updated all-at-once in this manner until either a time limit is reached or the improvement in the optimal value of \eqref{eq:slow_evaluation} is too small (or negative).  (In terms of evaluating the objective value of \eqref{eq:slow_evaluation}, it is important to evaluate the objective value explicitly, rather than consider the objective value reported by Ipopt, since the latter is influenced by internal relaxations of the bound constraints.)

For the first time that a contingency is evaluated, the active-set predictions for the complementarity constraints are initialized differently depending on the type of complementarity and contingency. With respect to active power (i.e., \eqref{eq.comp_active}) in a generator contingency, a one-dimensional bisection search is performed over the perturbation variable $\Delta_k$.  For each generator, a new production level is computed in a way that satisfies the complementarity constraints \eqref{eq.comp_active}.  The search is terminated when the overall added active power generation equals $1.01$ times the amount of active power lost by the generator contingency. The factor $1.01$ is used as an estimate that we would incur a 1\% increase in losses due to the rerouting of power flows. The active constraints resulting from this search are used to initialize the active-set predictions for each generator.  In all other cases, the active-set predictions are initialized to the middle segment.

We also tried reformulating the complementary constraints as penalty terms in the objective function.  Observe that each complementarity constraint in \eqref{eq.comp_active}--\eqref{eq.comp_reactive} can be written as $0 \leq \eta \perp \mu \geq 0$, which is equivalent to saying that $(\eta,\mu) \geq 0$ while $\eta\mu = 0$.  An approach that has been explored in the literature is to include the bound constraints $(\eta,\mu) \geq 0$ and add a term of the form $-\beta (\eta\mu)$ to the objective function, where $\beta$ is a positive penalty parameter \citep{HuRalph2004}.  This reformulation did not work well in our experiments, which we attribute to the negative curvature that such a penalty term introduces into the objective function.

\subsection{Fast Contingency Evaluation}\label{sec.fast}







Evaluating each contingency as described in the preceding subsections can be expensive computationally, especially due to the combinatorial nature of problem~\eqref{eq:slow_evaluation} from the presence of the complementarity constraints.  In an attempt to mitigate the computational costs that would be incurred by performing a full evaluation of each contingency during each iteration of the complementarity update loop described in the previous section, we employ a \emph{fast} contingency evaluation scheme that is able to quickly produce an upper bound on the optimal value of~\eqref{eq:slow_evaluation}.  This is done by solving a reduced contingency evaluation problem, then setting values of the remaining variables in a way that results in a feasible solution of problem~\eqref{eq:slow_evaluation}, which gives such an upper bound.

Each reduced problem involves the power flow equations~\eqref{eq.line_flow}--\eqref{eq.bus_balance}, but without the corresponding inequality constraints (i.e., physical bounds) and without the slack variables.  After fixing the controllable shunt susceptances (where, for our purposes, we use their values from the base case solution) and explicitly choosing the active complementarity segments, one obtains a square system of nonlinear equations for contingency $k \in \mathcal{K}$, which involves the equations:
\bequation \label{eq:fast_evaluation}
  \baligned
     \left\{
    \baligned
      a_k(u_k,y_k) &= 0 \\
      b_k(u_k,y_k) &= 0 \\
      \tau_k(u_0,y_0,u_k,y_k) &= 0.
    \ealigned
    \right.
  \ealigned
\eequation
Here, $\tau_k$ is derived from either $\tau_k^L$ or $\tau_k^R$ depending on the active-set prediction.  To solve \eqref{eq:fast_evaluation}, we use a Newton method; in particular, we call Ipopt to solve the system of equations.

After \eqref{eq:fast_evaluation} is solved, we update the active-set predictions for the complementarities based on any violated bounds, then project the voltage and generation values according to their bounds, adjust the remaining state variables, and compute slack variables accordingly to get a feasible solution to the full contingency problem.  One could stop there and employ the upper bound offered by this feasible solution.  In our approach, however, we continue the fast evaluation in an iterative manner.  In particular, using the updated active-set predictions for the complementarities, we repeat the aforementioned steps until the penalty value no longer decreases, a time limit is reached, or the penalty value falls below a cutoff value that is tightened iteratively by the overall algorithm.

If, for a given contingency, the penalty value obtained after termination of this fast evaluation process is greater than the cutoff value mentioned above, then we transition to the full evaluation of this contingency by solving~\eqref{eq:slow_evaluation} directly.  Otherwise, the contingency is deemed not to be worth a full evaluation at this stage in the overall solution algorithm.


\subsection{Avoiding Selection of Dominated Contingencies}\label{sec.dominated}




Iteratively throughout our solution algorithm, we select a specified number (three, in our implementation) of additional contingencies to include in subsequent master problems.  While our choice of contingencies is informed by their associated penalty values, simply selecting the contingencies with the largest penalties may result in redundant contingencies being added to the master problem.  In other words, adding the constraints associated with one contingency may implicitly result in the (near) satisfaction of all constraints associated with another contingency (or contingencies).  Hence, the latter contingencies are \emph{dominated} by the former.

Contingencies consisting of failures of nearby generators, transmission lines, or transformers can often be (nearly) redundant.  Some of these redundant contingencies can be identified easily.  For instance, the failure of an identical parallel line or identical generator at the same bus is clearly redundant since the failure of one unit will have the same impact on the system as the failure of its counterpart. 
%
In our solution algorithm, we preprocess the set of contingencies to eliminate these trivially redundant contingencies before we start our algorithm.

Identifying other near-redundancy in contingencies is more challenging.  However, since each contingency added to the master problem introduces a significant computational burden, it is worthwhile to make efforts to avoid inclusion of nearly redundant contingencies.  To this end, \citep{capitanescu2007contingency} introduced the concept of an \emph{individually dominated contingency}, which attempts to identify a redundant contingency by observing its corresponding constraint violations.


\begin{definition}\citep{capitanescu2007contingency}
  \textbf{Individually Dominated Contingency.}\\
  \emph{
  Contingency $k \in \Kcal$ is \textbf{individually dominated} by contingency $j \in K$ if
  \bequation \label{eq:IDC}
    \baligned
      \varsigma_{j} \geq \varsigma_{k}, \quad \varsigma_{j}^+ \geq \varsigma_{k}^+, \quad \varsigma_{j}^- \geq \varsigma_{k}^- \quad \text{(component-wise)}
    \ealigned
  \eequation 
  for the optimal penalty variables.}
\end{definition}
The approach in~\citep{capitanescu2007contingency} enforces the constraints associated with all contingencies that are not individually dominated. While somewhat effective, our initial numerical experiments suggested that this approach is too conservative for the SC-AC-OPF problem considered in the GO Challenge 1 competition since too many nearly redundant contingencies (i.e., with nonzero, but small penalty values) are still added to the master problem. Therefore, we relax condition~\eqref{eq:IDC} by only considering the most violated constraint for each contingency, rather than all constraints, through the introduction of the concept of a \emph{maximum violation dominated contingency}.

\begin{definition}
  \textbf{Maximum Violation Dominated Contingency.}\\
  \emph{
  Let $\overline{i}_k$ be the index of the constraint in contingency $k \in \Kcal$ with the largest violation (i.e., the index of the constraint with the largest associated slack variable $\varsigma_k$, $\varsigma_k^+$, or $\varsigma_k^-$), and let $\bar{\varsigma}_k$ be the corresponding violation. We say contingency $k$ is \textbf{maximum violation dominated} by contingency $j \in \Kcal$ if
  \bequation \label{eq:CDC}
    \baligned
      \overline{i}_j = \overline{i}_k \quad \text{and} \quad  \bar{\varsigma}_j > \bar{\varsigma}_k.
    \ealigned
  \eequation
  }
\end{definition}

During the selection of contingencies to include in subsequent master problems, we make sure to select a set of contingencies that are not maximum violation dominated.

\subsection{Preprocessing and Handling Degeneracy}\label{sec.preprocessing}

To achieve faster convergence for the interior-point method, we preprocess the input data to eliminate degeneracies. Some of these degeneracies result from parallel lines that have identical electrical parameters (i.e., $g_e$, $b_e$, and $\overline{R}_e$). Since they have the same terminal voltages, the power flows on these parallel lines are also identical. Hence, the line flow limits \eqref{eq.line_current_rating} for these lines are redundant and all but one can be eliminated. Flow limits for identical parallel transformers are similarly redundant. However, when dealing with the contingency associated with one of these lines, we have to make sure that at least one copy of the flow limit constraints is present when one of the lines is removed from the network. 



Our preprocessing step also eliminates degeneracies associated with generators that are located at the same bus. In particular, we reduce the number of variables and variable bounds by aggregating the reactive power outputs of all generators at the same bus.


\subsection{Overall Algorithm}\label{sec.overall}

Our overall solution algorithm consists of three major components, all occurring in concert in parallel: solving master problems, evaluating contingencies, and applying ranking schemes to identify important contingencies.  
The algorithm maintains a priority list of contingencies that determines in which order they should be evaluated and which top three should be included in subsequent master problems.  This list is updated at various steps of our algorithm.

Our algorithm is summarized in the steps below, which we refer to as \ouralg{}.  It is implemented as \emph{Code 1}, which is depicted on the left-hand side of Figure~\ref{fig.flow_chart}.

\begin{enumerate}
  \item Remove trivially redundant contingencies (Section~\ref{sec.dominated}) and preprocess the data (Section~\ref{sec.preprocessing}).
  \item Solve the base case problem to obtain a base case solution.  For this problem, the variables are initialized with a \emph{flat start}: active power generation values are set to their upper limits; reactive power generation values are set to zero; voltages and controllable shunt susceptances at the buses are initialized to the midpoints within their described limits; and line and transformer flow variables and voltage angles are initialized to zero.
  \item Write the current base case solution to a file.
  \item\label{alg:first_eval} Initialize the contingency priority list (Section~\ref{sec.initial_ranking}).
  \item\label{alg:init_fast} Apply fast contingency evaluation (Section \ref{sec.fast}) to the contingencies, including complementarity updates (Section~\ref{sec.complementarities}), according to the priority list until a time limit (of one minute) is reached or until all contingencies have been evaluated. Re-sort the priority list in order from largest to smallest penalty value.
  \item \label{alg:beginnig_of_loop} Perform full contingency evaluations (Section~\ref{sec.contingency_evaluation}) starting from the top of the list, given an overall time limit of 30 seconds.  Re-sort the priority list in order from largest to smallest penalty value.
  \item\label{s:addcont} Add the top three contingencies from the list to the master problem, avoiding dominated contingencies (Section~\ref{sec.dominated}).
  \item Solve the master problem, with no adjustments to the complementarity segments, while at the same time continuing to evaluate contingencies in the order of the priority list, with preference given to those that have not yet been evaluated.  Each evaluation of a contingency starts with the fast evaluation (Section~\ref{sec.fast}) including complementarity updates.  If the penalty value after fast evaluation is still above a threshold, then full evaluation (Section~\ref{sec.contingency_evaluation}) with complementarity updates (Section~\ref{sec.complementarities}) is performed.
  \item On completion of the solve of the master problem, write the new base case solution to a file and terminate the concurrent contingency evaluation procedure.
  \item Update the priority list based on the most recently computed penalty values.  If the priority list is not empty and time remains, return to Step~\ref{alg:beginnig_of_loop}.  Otherwise, stop.

\end{enumerate}

As previously mentioned, in the GO Challenge~1 competition, \emph{Code 1} was given a time limit of 10 minutes or 45 minutes, depending on the competition category, and terminated once the time limit was exceeded.  The most recent base case solution, written to a file, is the product of \emph{Code~1}.  
Our implementation of \emph{Code~2} is depicted in the right-hand side of Figure~\ref{fig.flow_chart}.  After the base case solution produced by \emph{Code~1} is read from the file, all contingencies are evaluated in parallel.  For each contingency, the fast evaluation in Section~\ref{sec.fast} including complementarity updates is used first.  If the penalty value after this evaluation is above a threshold, it is further reduced by a full contingency evaluation as described in Section \ref{sec.contingency_evaluation}.

For the competition, it was crucial to compute results for \emph{all} contingencies within the time limit for \emph{Code~2}, set by the competition organizers to (2 $\times$ number of contingencies) seconds.  To ensure timely termination, at the beginning of a new contingency evaluation, we calculate how much total computing time is left, accumulated over all threads combined, and divide this total remaining time by the number of unevaluated contingencies. If from then on all contingencies would require that much time, it would be guaranteed that \emph{Code 2} would finish within the overall time limit for \emph{Code 2}.  However, since most contingency evaluations are finished faster, the remaining time per contingency is expected to increase over the run of \emph{Code 2}, with the largest time limits observed at the very end. To make best use of this strategy, we process the contingencies in the reverse order obtained by the ranking procedure described in Section~\ref{sec.initial_ranking}. Since contingencies for which the evaluation requires more time are typically ranked earlier in our list, they are processed at the final stages of \emph{Code 2} and thus receive the largest time limit.

\section{Numerical Results}
\label{sec.results}

In this section, we present the results of experiments to demonstrate the effectiveness of a few of our algorithmic strategies proposed in Section~\ref{sec.algorithm}.  Our \emph{Code~1} and \emph{Code~2} have been implemented in C++ with MPI; recall Figure~\ref{fig.flow_chart} for a visualization of our solution algorithm.

The nonlinear interior-point optimization software known as Ipopt \citep{WaecBieg06} was used to solve the base case problem and the master problems with selected contingencies; these problems all have the form in \eqref{prob.high_level}.  Ipopt was also used to solve the contingency evaluation problems as described in Section~\ref{sec.contingency_evaluation}, as well as to solve the power flow equations for the fast contingency evaluation strategy described in Section~\ref{sec.fast}.
During the competition, we used an asynchronous parent--child paradigm with hierarchical parallelism.  Communication between the six processes, each one executed on a different compute node, was done via MPI.  Twenty pthreads were run within each process.
The parent process ran only the base case and master problem optimizations even though these required only one thread.  Since these are the most crucial tasks in terms of the time limit, reserving an entire compute node to this thread meant that it was not slowed down by sharing hardware resources with other threads.  The child processes were responsible for the contingency evaluations.

For the experiments presented in this section, we did not have access to the GO Competition computing platform.  Instead, we ran our codes on a single 
Linux server node with two Intel Xeon CPUs running 40 (hyperthreaded) threads each and 256GB RAM.

Our algorithm involves several heuristic strategies.  The experiments here are intended to demonstrate the contribution of each of the following of these heuristics: 
initial contingency ranking (see Section~\ref{sec.num_initial_ranking}), complementarity update (see Section~\ref{sec.num_complementarities}), fast evaluation (see Section~\ref{sec.num_fast}), and contingency selection (see Section~\ref{sec.num_cont_selection}), with different heuristic choices being switched on and off. 

The test cases in our experiments were taken from the Challenge 1 offline datasets \citep{GoData}. We evaluated all scenarios for the selected networks. 
No time limits were imposed during these experiments because the hardware used was significantly less powerful than that available during the competition. In addition, since our purposes here are to isolate the effect of each of the aforementioned heuristics, we want to show results that are unaffected by a time limit.


\subsection{Comparison: Contingency Initial Rankings}\label{sec.num_initial_ranking}

We employ a regression model to provide an initial ranking of contingencies; recall Section~\ref{sec.initial_ranking}. This ranking prescribes the order in which the contingencies are initially evaluated in Step~\ref{alg:first_eval} of \ouralg{}. In the competition, only a limited number of contingencies could be evaluated within the time constraint. Thus, it is critical that the most important contingencies are near the top of the initial contingency ranking so that they are the ones that are evaluated within the time constraints. 
In this section, we compare the heuristic we included in the algorithm with three other simpler heuristics that involve fewer features.  In particular, we compare the following four ranking heuristics using the initial base case solution:
\begin{itemize}
    \item Heuristic $l^p$: rank the contingencies by the associated amount of active power loss ($p_g$ for a generator contingency and $\max\{p^o_e, p^d_e\}$ for a line/transformer contingency).
    \item Heuristic $l^s$: rank the contingencies by the associated amount of apparent power loss ($\sqrt{p_g^2 + q_g^2}$ for a generator contingency and $\max\{\sqrt{(p^o_e)^2 + (q^o_e)^2}, \sqrt{(p^d_e)^2+(q^d_e)^2}\}$ for a line (and similarly for a transformer) contingency).
    \item Heuristic $l^c$: rank the contingencies by the associated apparent power loss relative to the corresponding capacity \Big($\sqrt{\frac{p_g^2 + q_g^2}{\overline{p}_g^2 + \overline{q}_g^2}}$ for a generator contingency and $\max\left\lbrace \sqrt{\frac{(p^o_e)^2 + (q^o_e)^2}{\overline{R}_e v_{i^o_e}}}, \sqrt{\frac{(p^d_e)^2+(q^d_e)^2}{{\overline{R}_e v_{i^d_e}}}}\right\rbrace$ for a line (and similarly for a transformer) contingency\Big).
    \item Heuristic $l^w$: rank the contingencies by the regression model as in Section \ref{sec.initial_ranking}.
\end{itemize}

To compare the effectiveness of these alternatives, we solved the base case problem and obtained the initial rankings with each of these ranking schemes.  As a base line, we also performed full penalty minimization by solving \eqref{eq:slow_evaluation} for each contingency to obtain their smallest penalty values. 

\begin{figure}
  \begin{subfigure}[t]{.45\textwidth}
    \centering
    \includegraphics[width=\linewidth]{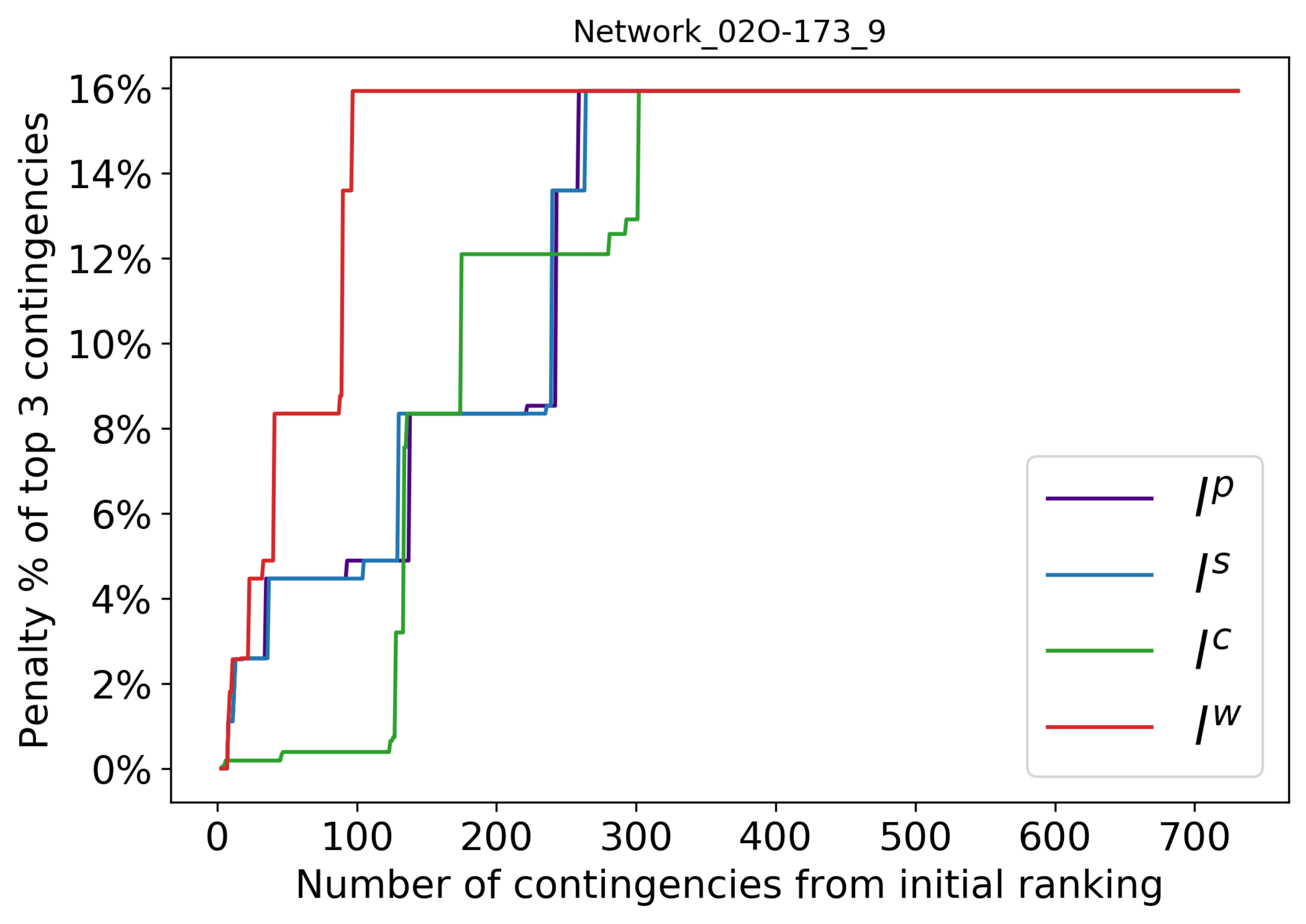}
    \caption{Network\_02O-173\_9}
  \end{subfigure}
  \hfill
  \begin{subfigure}[t]{.45\textwidth}
    \centering
    \includegraphics[width=\linewidth]{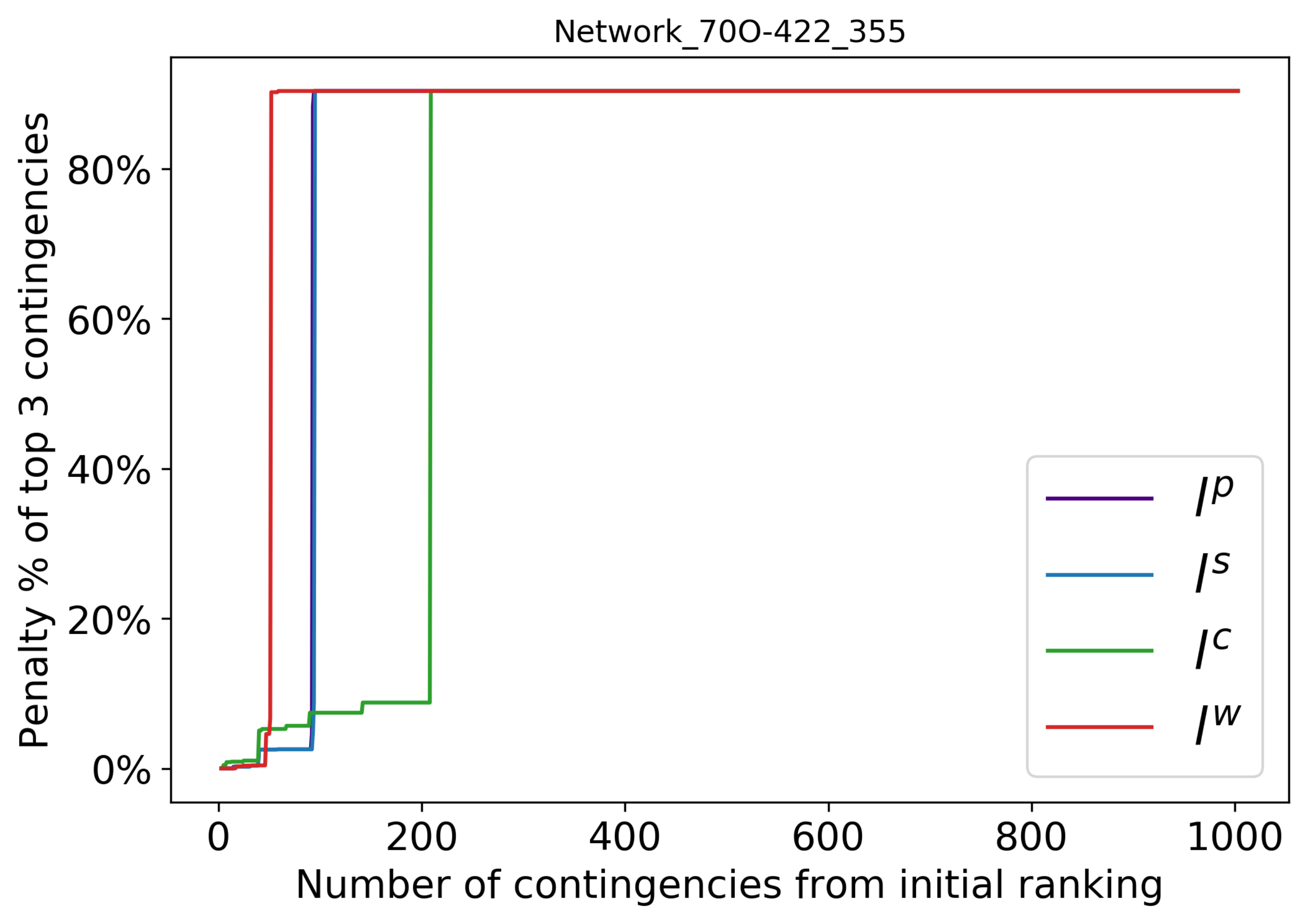}
    \caption{Network\_70O-422\_355}
  \end{subfigure}

  \medskip

  \begin{subfigure}[t]{.45\textwidth}
    \centering
    \includegraphics[width=\linewidth]{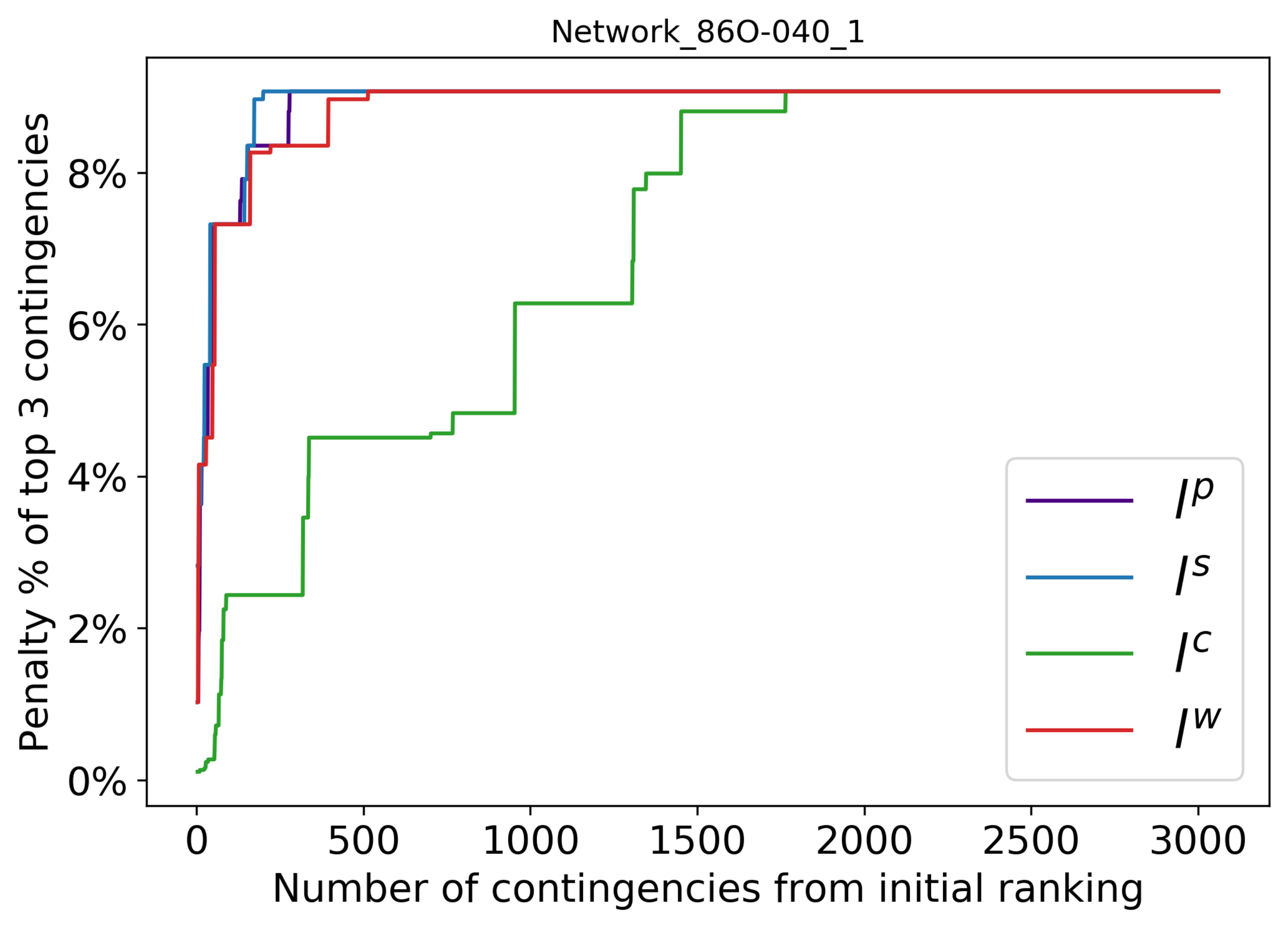}
    \caption{Network\_86O-040\_1}
  \end{subfigure}
  \hfill
  \begin{subfigure}[t]{.45\textwidth}
    \centering
    \includegraphics[width=\linewidth]{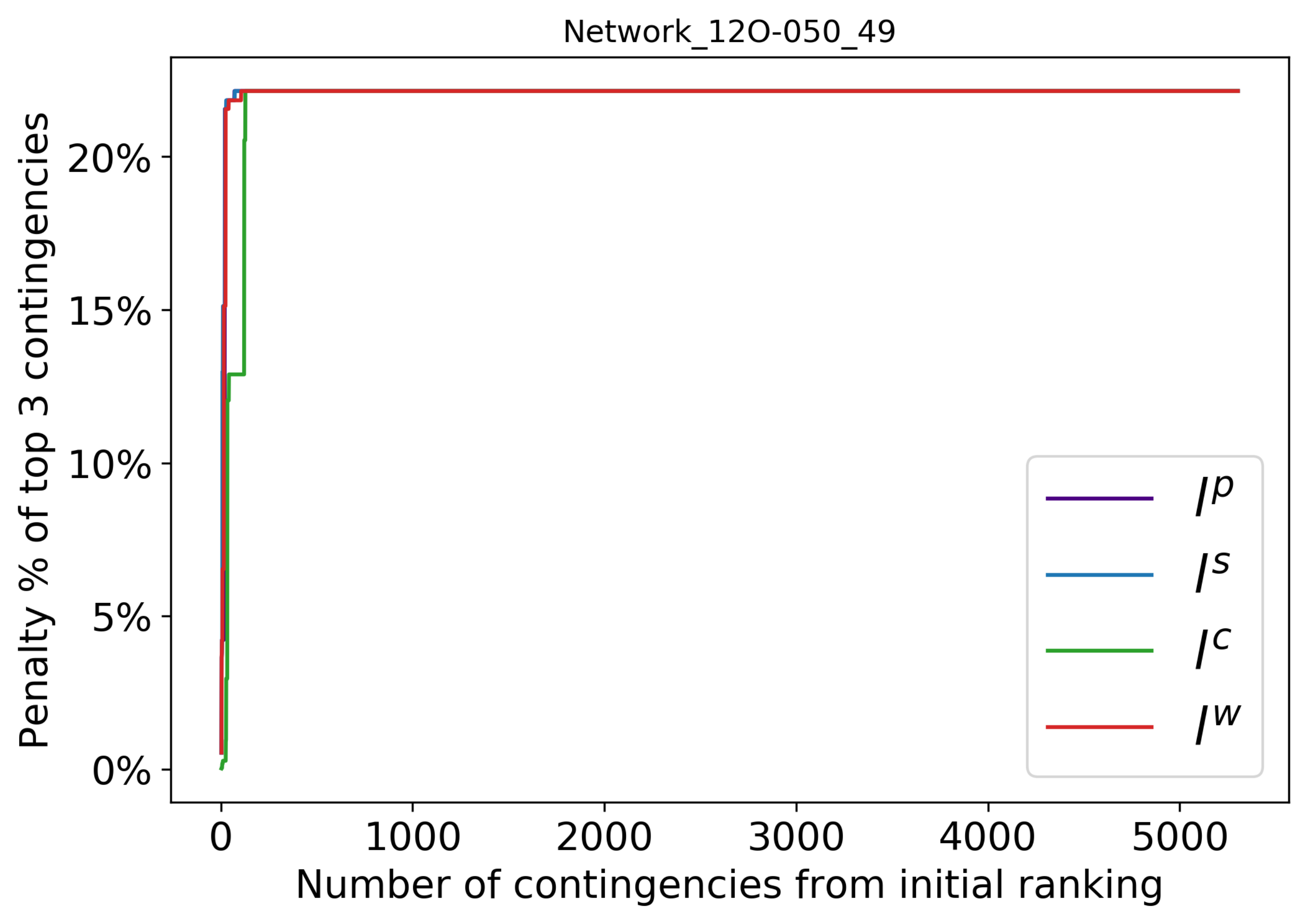}
    \caption{Network\_12O-050\_49}
  \end{subfigure}
  \caption{Comparison of initial contingency ranking schemes.}
  \label{fig.init_ranking}
\end{figure}

In Step~\ref{s:addcont} of \ouralg{}, three contingencies are selected after the contingency evaluation. Therefore, we only focus on the predictions of the top three contingencies using the above four ranking schemes. From the top of the ranking generated by each heuristic, we compute the penalty values of the top three contingencies as a percentage of the total penalty values of all contingencies. Figure~\ref{fig.init_ranking} compares the ranking heuristics for four networks with different sizes. One representative scenario is selected for each network.  On the horizontal axis, we increase the number of contingencies that are considered from the initial ranking.  For each number of contingencies considered, we determine which are the top three using the different ranking schemes and compute the aforementioned percentage.  In this manner, the penalty percentage monotonically increases as one moves to the right in the horizontal axis.  A jump in the graph for a given heuristic indicates that a contingency with large penalty value has been identified once one more contingency is considered.

Figure~\ref{fig.init_ranking} shows that all four ranking schemes identify important contingencies at an early stage of the contingency evaluation.  However, our ranking heuristic $l^w$ generally captures more important contingencies earlier in the list, as these representative figures illustrate.  Overall, given any number of contingencies in the ranked list that have been evaluated, the top three contingencies identified via the regression model consistently have higher penalty values compared to those identified by the other three heuristics.

In the competition, the parallel execution of the initial penalty evaluations in Steps \ref{alg:init_fast} and \ref{alg:beginnig_of_loop} made it possible to compute the correct (minimal) penalty value of several hundreds of contingencies within the time limits.  The results in Figure~\ref{fig.init_ranking} indicate that the most important contingencies could often be correctly identified and were added to the first master problem.

\subsection{Comparison: Complementarity Updates}\label{sec.num_complementarities}
The choices of the active segments for the complementarity constraints \eqref{eq.comp_active} and \eqref{eq.comp_reactive} are updated using an active-set approach during contingency evaluation; recall Section~\ref{sec.complementarities}. To illustrate the effectiveness of this approach, we next compare the contingency evaluation results obtained with and without complementarity updates.

In this experiment, we solve the base case problem (without contingencies) using \emph{Code 1}, then perform contingency evaluations using \emph{Code 2} with the base case solution from \emph{Code 1}.  In this contingency evaluation, we compare two approaches:
\begin{enumerate}
    \item Without complementarity updates: each contingency is evaluated by solving \eqref{eq:slow_evaluation} using the default complementarity initialization.
    \item With complementarity updates: adjustments to the complementarity constraints are made iteratively within the contingency evaluations until the decrease in the penalty value is sufficiently small (as described in Section~\ref{sec.complementarities}).
\end{enumerate}

\bfigure 
  \centering
  \includegraphics[width=0.7\textwidth]{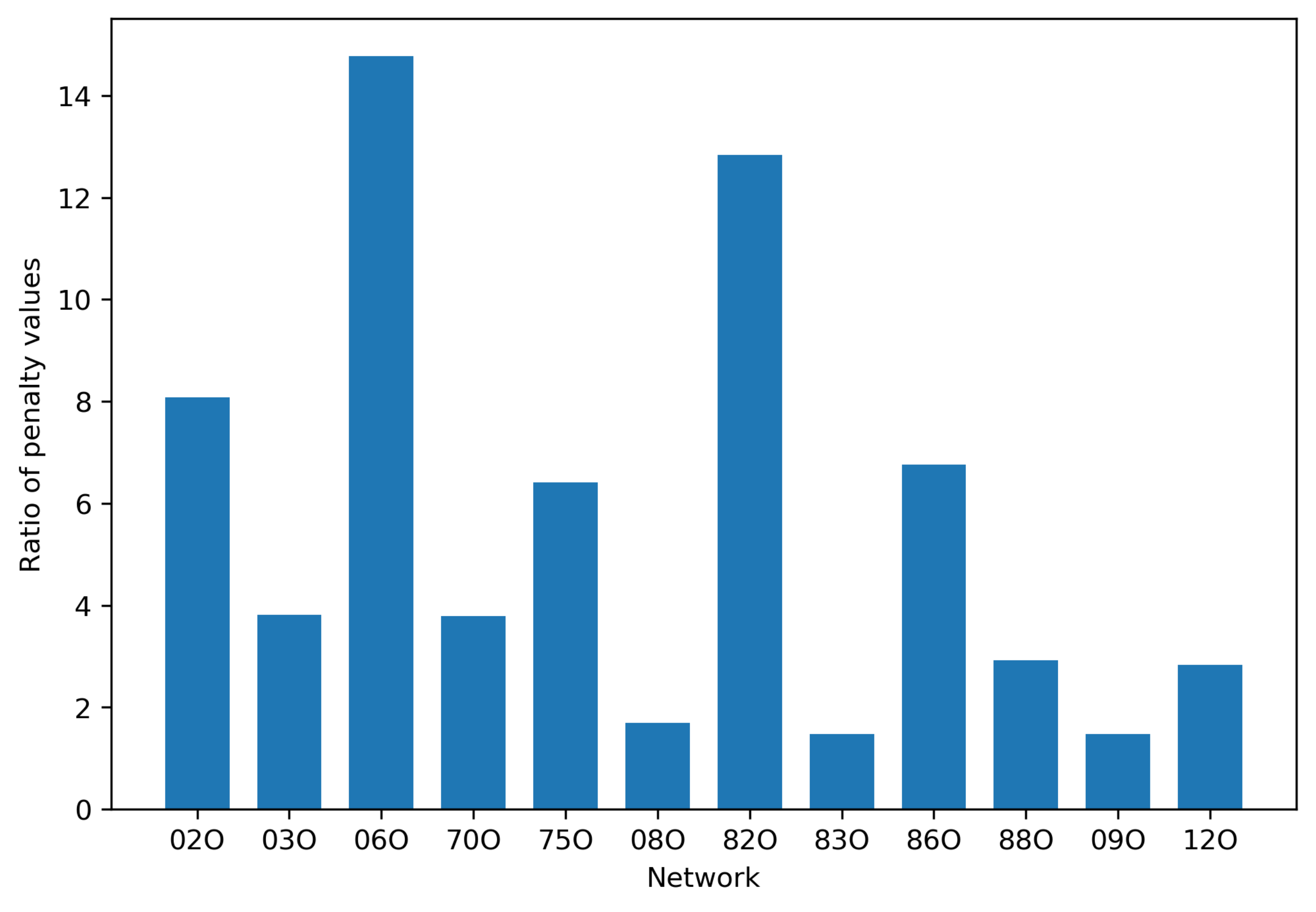}
  \caption{Ratio of penalty values obtained without complementarity updates to penalty values obtained with complementarity updates.}
  \label{fig.comp_comparsion}
\efigure

Figure~\ref{fig.comp_comparsion} shows the factors by which the penalty values obtained without complementarity updates are worse than those obtained with complementarity updates. The result for each network is obtained by averaging over the scenarios.  Note that since the complementarity updates allow for further optimization in the contingency evaluation, the penalty values obtained with complementarity updates are always less than the penalty values obtained without complementarity updates; hence, each bar in the figure is greater than or equal to 1.  Allowing complementarity updates often results in significantly reduced penalty values for each network.

\subsection{Comparison: Fast Contingency Evaluation}\label{sec.num_fast}

Our fast contingency evaluation scheme (see Section~\ref{sec.fast}) is intended to mitigate the heavy computational cost required for full contingency evaluation.  Indeed, in many cases, fast evaluation serves as a pre-screening step for filtering out feasible contingencies (i.e., with low penalty) from requiring a full evaluation.  To assess the computational improvements resulting from this pre-screening step, we conduct a comparison of contingency evaluation with and without fast evaluation.  Similar to the preceding experiment, we perform contingency evaluation in \emph{Code 2} based on the base case solution obtained from \emph{Code 1}. Two approaches are compared here:
\begin{enumerate}
    \item With pre-screening: fast evaluation is performed first, followed by full evaluation if the penalty value obtained after fast evaluation is larger than a cutoff value.
    \item Without pre-screening: each contingency is fully evaluated without fast evaluation.
\end{enumerate}
We set the cutoff value for both fast and full evaluations such that the constraint violations are less than $2\times 10^{{-}2}$ per unit. (This constraint violation corresponds to the first breakpoint in the piecewise linear penalty function used in the competition.) 
Note that complementarity updates are included in both approaches. 

\bfigure 
  \centering
  \includegraphics[width=0.7\textwidth]{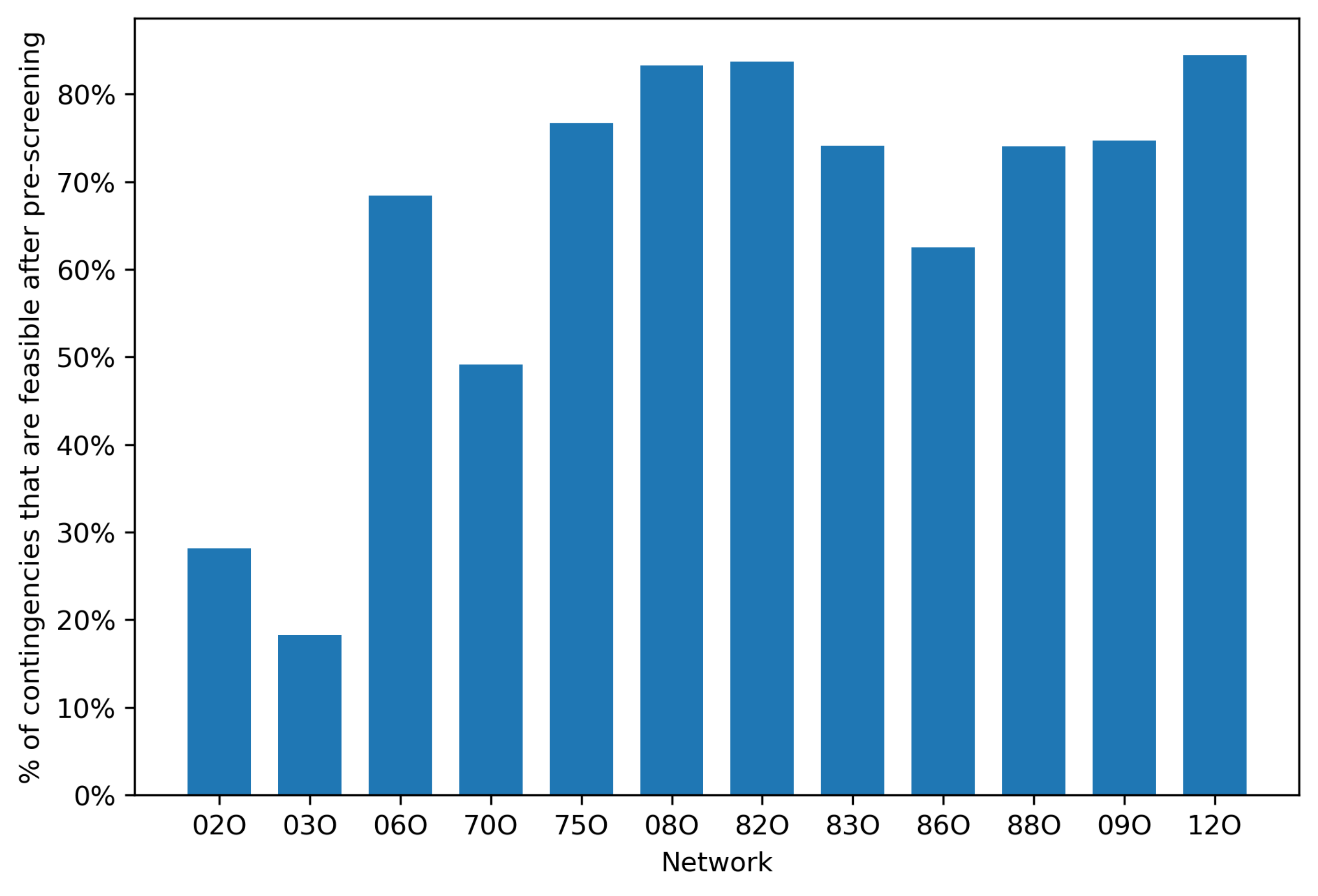}
  \caption{Percentage of contingencies with penalty value below a cutoff after pre-screening.}
  \label{fig.fasteval_num}
\efigure

\bfigure 
  \centering
  \includegraphics[width=0.7\textwidth]{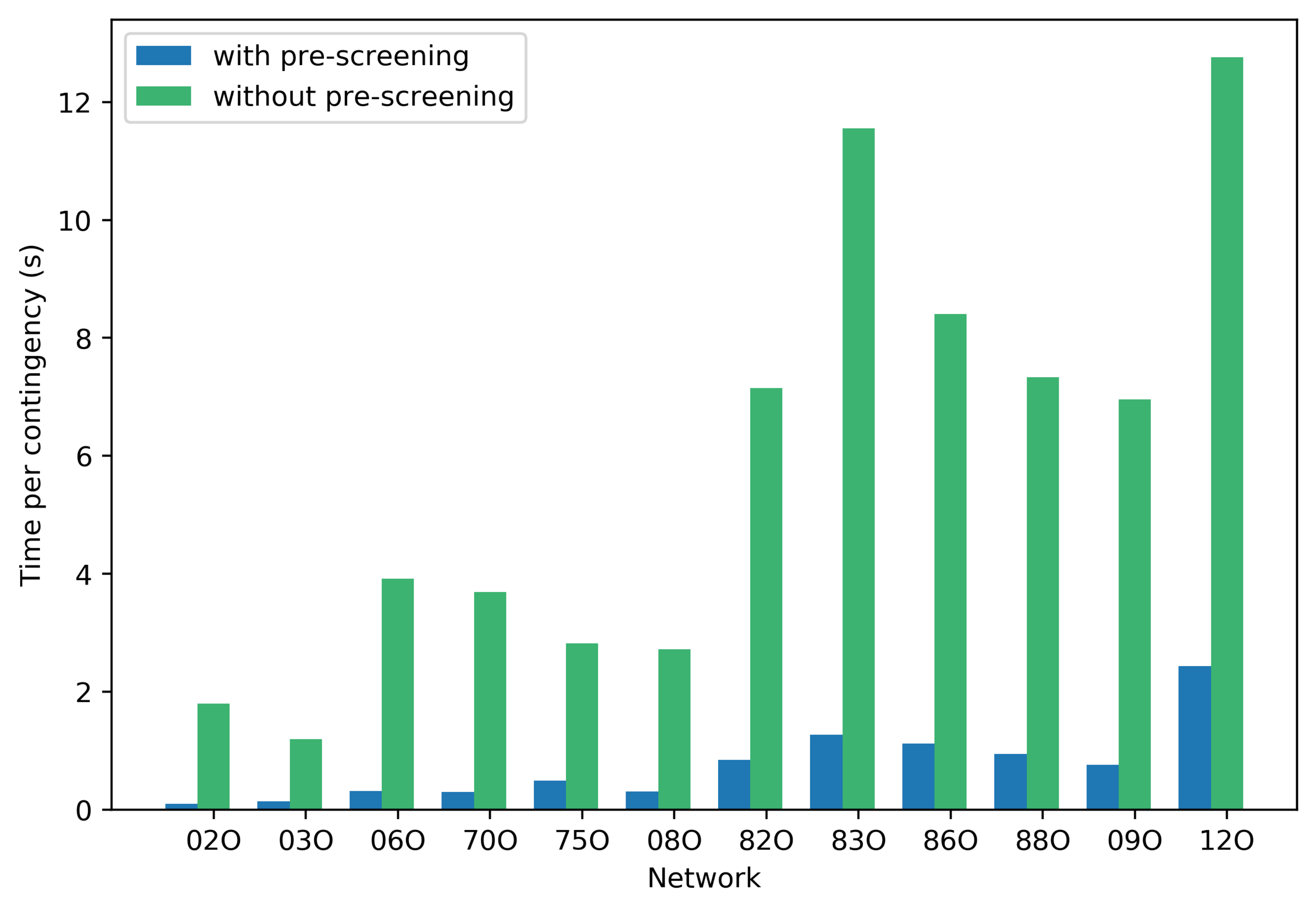}
  \caption{Comparison of average evaluation time among contingencies that are determined to be (sufficiently) feasible after pre-screening.}
  \label{fig.fasteval_comparsion_time}
\efigure

\bfigure 
  \centering
  \includegraphics[width=0.7\textwidth]{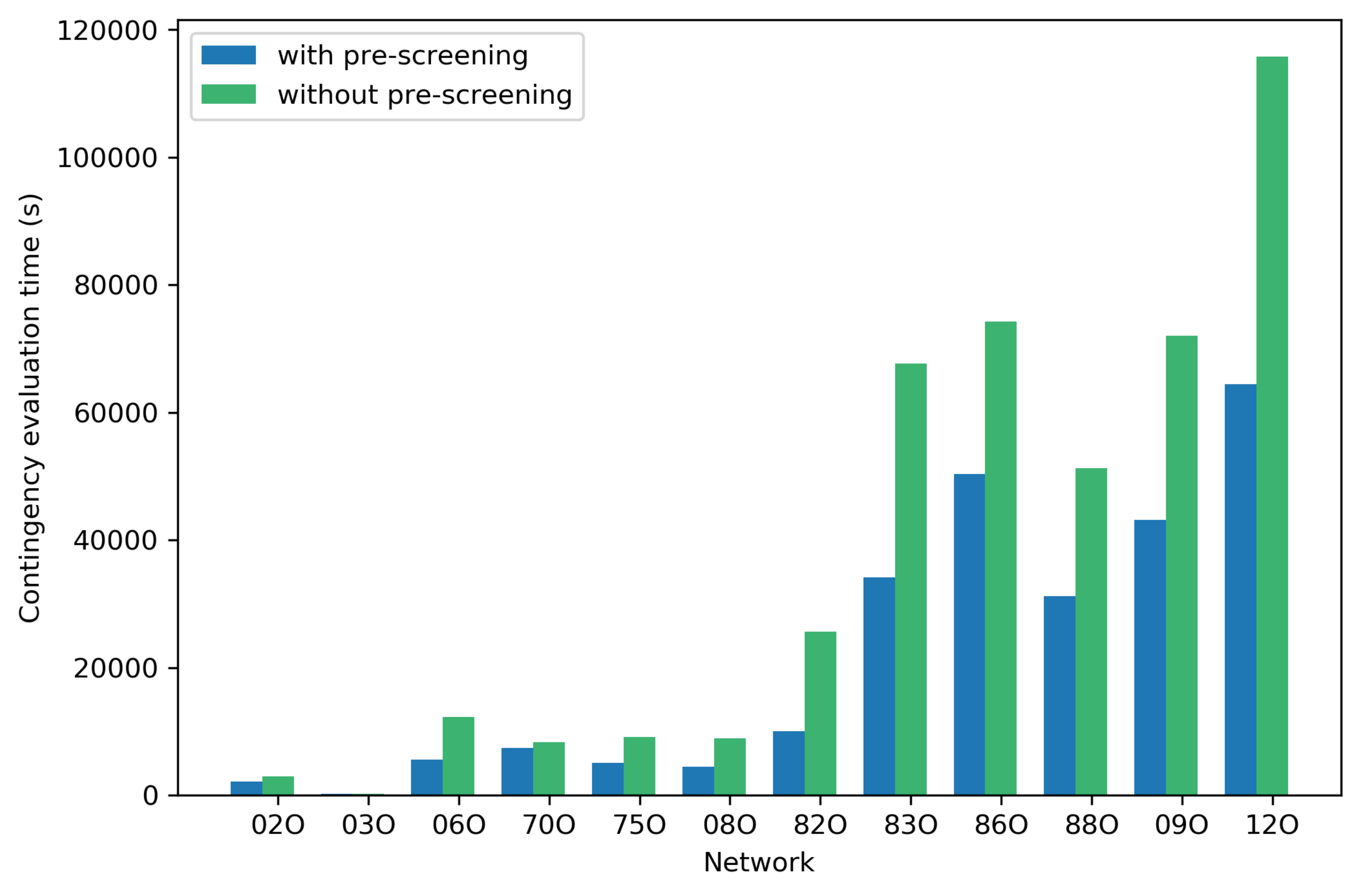}
  \caption{Comparison of contingency evaluation times with and without pre-screening.}
  \label{fig.fasteval_comparsion_total}
\efigure

Figure~\ref{fig.fasteval_num} shows the percentage of contingencies that are feasible after pre-screening. For most of the networks, fast evaluation eliminates over half of the contingencies.
Among those contingencies that are feasible after pre-screening, we compare the average computation times with and without pre-screening in Figure~\ref{fig.fasteval_comparsion_time}. This figure shows significant improvement in computational time through pre-screening. 
Figure~\ref{fig.fasteval_comparsion_total} compares the total contingency evaluation times with and without pre-screening. This figure shows that pre-screening allows the evaluation of more contingencies within a limited time, thus identifying more important contingencies to include in the master problem.


\subsection{Comparison: Contingency Selection}\label{sec.num_cont_selection}

As discussed in Section~\ref{sec.algorithm}, it is important to identity the most important contingencies to include in the master problem. Even though the initial ranking scheme gives a good approximation of the true penalties, the top contingencies captured by the initial ranking might not be the most important. Therefore, it is necessary to evaluate the contingencies to identify those that are most important.  In \ouralg{}, fast evaluation is applied first to contingencies in the order of the initial ranking.  This eliminates a large percentage of unimportant contingencies within a short amount of time. The algorithm proceeds with full evaluation only for the remaining contingencies. The dominance filtering technique proposed in Section~\ref{sec.dominated} helps avoid adding redundant contingencies into the master problem.

To show how each step in this procedure helps in identifying important contingencies, we compare four different contingency selection approaches, each of which adds an additional step to the previous approach:

\begin{enumerate}
    \item Initial ranking: contingencies are ranked using our initial ranking scheme.
    \item Fast evaluation: fast evaluation is applied to the contingencies in the order of the initial ranking. Contingencies are then ranked by their penalty values computed by fast evaluation.
    \item Full evaluation: full evaluation is performed if the penalty value resulting from fast evaluation is above a cutoff value. Contingencies are then ranked by their penalty values computed by full evaluation.
    \item Dominance filtering: the dominance filtering technique is applied after the full evaluation of the contingencies. Dominated contingencies are eliminated from the ranking.
\end{enumerate}

For each of these approaches, we augment the base case with the constraints from the top three contingencies to form four master problems. We then solve these master problems and evaluate all contingencies for the new base case solution to compute the corresponding objective values including the penalty costs. To characterize the quality of the solutions obtained by just including three contingencies in the master problem, we compare these objective values with the best score from the final round of the competition. The results for each network are obtained by taking the average of the results over each scenario. Figure~\ref{fig.cont_selection} shows the optimality gaps between these scores and the best scores from the competition.  In nearly all cases, sequentially augmenting the contingency selection scheme to include the initial ranking strategy, fast evaluation, full evaluation, and dominance filtering monotonically reduces the total cost, in some cases significantly.


\bfigure 
  \centering
  \includegraphics[width=0.7\textwidth]{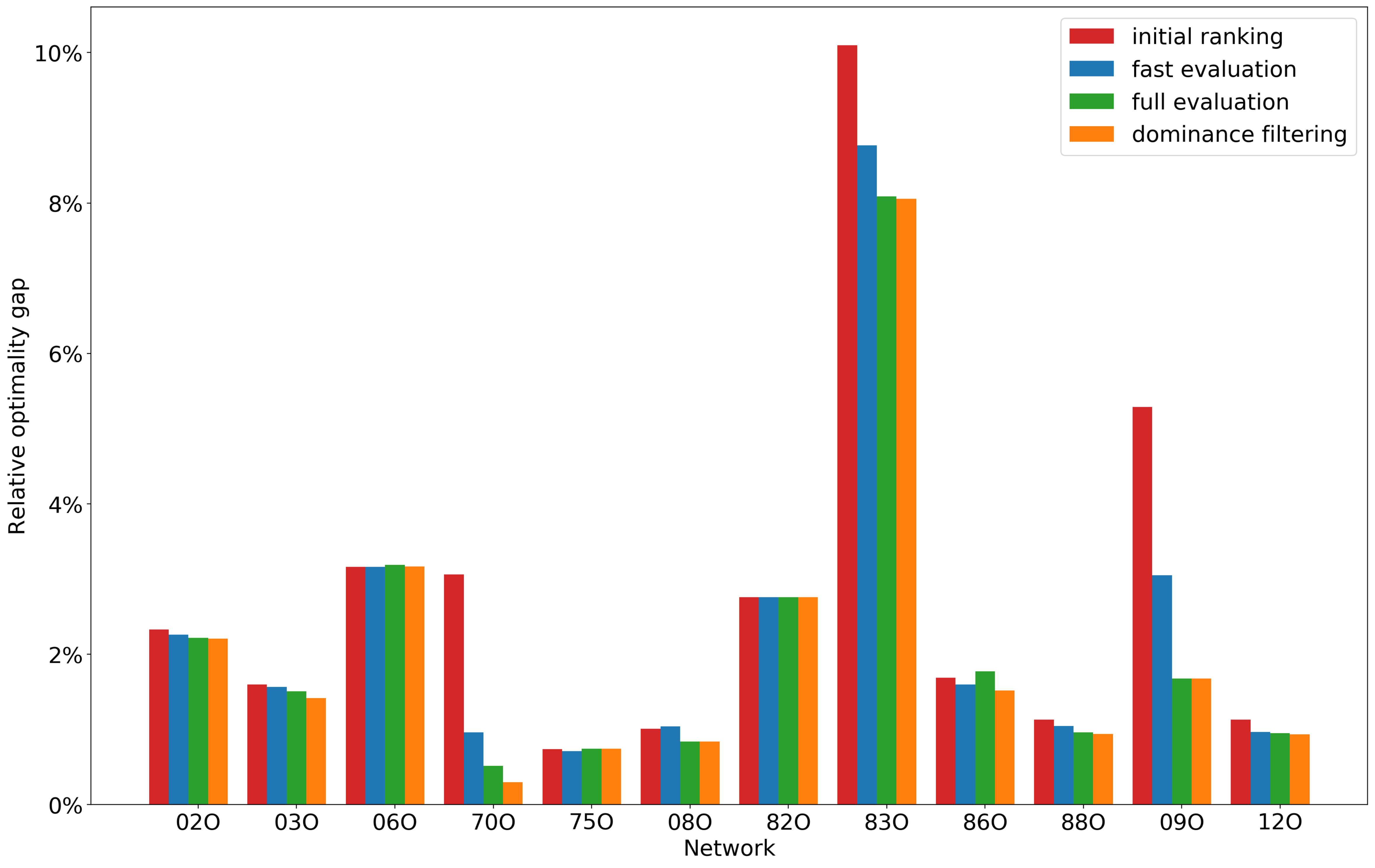}
  \caption{Comparison of combined contingency selection schemes.}
  \label{fig.cont_selection}
\efigure

\section{Conclusion}\label{sec.conclusion}

This paper presents the algorithm that our team, ``GO-SNIP,'' submitted for the ARPA-E GO Competition, Challenge 1.  The key features of our algorithm are strategies for contingency selection, fast contingency evaluation, handling complementarity constraints, avoiding issues related to degeneracy, and exploiting parallelism.  The results of the numerical experiments that are presented isolate the effectiveness of important features of our solution algorithm.  Our submission for the GO Challenge~1 competition, along with those of the other participating teams, demonstrate that tailored nonlinear optimization techniques are capable of solving large-scale SC~AC-OPF problems within industrially relevant time limits.

Ongoing efforts to improve upon the SC~AC-OPF algorithms developed for the ARPA-E Grid Optimization Challenge~1 competition benefit from recently developed extensive repositories of large-scale realistic power system datasets. References to the specific datasets used in Challenge~1 and the currently ongoing Challenge~2 of the GO Competition are provided in \citep{GoData} and \citep{GoData2}, respectively. Other large repositories of power system data include DR~POWER \citep{drpower} and The GRID~DATA Repository \citep{bettergrids}.  Additionally, researchers often use the curated set of AC-OPF problems in the PGLib-OPF repository for algorithmic benchmarking purposes~\citep{pglib}.